\documentclass[12pt]{article}
\usepackage{amsmath, amsthm, amssymb, fullpage}
\usepackage[all]{xy}

\newtheorem{theorem}{Theorem}[subsection]
\numberwithin{equation}{theorem}
\newtheorem{lemma}[theorem]{Lemma}
\newtheorem{cor}[theorem]{Corollary}
\newtheorem{prop}[theorem]{Proposition}
\theoremstyle{definition}
\newtheorem{defn}[theorem]{Definition}
\newtheorem{example}[theorem]{Example}
\newtheorem{remark}[theorem]{Remark}
\newtheorem{hypothesis}[theorem]{Hypothesis}
\newtheorem{notation}[theorem]{Notation}

\newcommand{\AAA}{\mathbb{A}}
\newcommand{\GG}{\mathbb{G}}
\newcommand{\PP}{\mathbb{P}}
\newcommand{\QQ}{\mathbb{Q}}
\newcommand{\RR}{\mathbb{R}}
\newcommand{\ZZ}{\mathbb{Z}}
\newcommand{\calE}{\mathcal{E}}
\newcommand{\calF}{\mathcal{F}}
\newcommand{\calO}{\mathcal{O}}
\newcommand{\calR}{\mathcal{R}}
\newcommand{\gotho}{\mathfrak{o}}

\newcommand{\del}{\partial}

\newcommand{\dual}{\vee}
\newcommand{\be}{\mathbf{e}}
\newcommand{\bv}{\mathbf{v}}

\DeclareMathOperator{\charac}{char}
\DeclareMathOperator{\dom}{dom}
\DeclareMathOperator{\epi}{epi}
\DeclareMathOperator{\Frac}{Frac}
\DeclareMathOperator{\Gal}{Gal}
\DeclareMathOperator{\GL}{GL}
\DeclareMathOperator{\intdom}{intdom}
\DeclareMathOperator{\inte}{int}
\DeclareMathOperator{\rank}{rank}
\DeclareMathOperator{\ratrank}{ratrank}
\DeclareMathOperator{\sep}{sep}
\DeclareMathOperator{\Spec}{Spec}
\DeclareMathOperator{\spect}{sp}
\DeclareMathOperator{\trdeg}{trdeg}
\DeclareMathOperator{\unr}{unr}

\begin{document}

\title{Semistable reduction for overconvergent $F$-isocrystals, III:
Local semistable reduction at monomial valuations}
\author{Kiran S. Kedlaya \\ Department of Mathematics, Room 2-165 \\ Massachusetts
Institute of Technology \\ 77 Massachusetts Avenue \\
Cambridge, MA 02139 \\
\texttt{kedlaya@mit.edu}}
\date{June 12, 2008}

\maketitle

\begin{abstract}
We resolve the local semistable reduction problem for 
overconvergent $F$-isocrystals
at monomial valuations (Abhyankar valuations of height 1 and
residue transcendence degree 0).
We first introduce a higher-dimensional analogue of the generic radius
of convergence for a $p$-adic differential module, which obeys a
convexity property. We then combine this convexity property with
a form of the $p$-adic local monodromy theorem for so-called fake annuli.
\end{abstract}

\tableofcontents

\section{Introduction}

This paper is the third of a series, preceded by \cite{kedlaya-part1, 
kedlaya-part2}.
The goal of the series is to prove a 
``semistable reduction'' theorem for overconvergent $F$-isocrystals,
a class of $p$-adic analytic objects associated to schemes of finite type
over a field of characteristic $p>0$.
Such a theorem is expected to have consequences for the theory of rigid
cohomology, in which overconvergent $F$-isocrystals play the role of
coefficient objects.

In \cite{kedlaya-part1}, it was shown that the problem of extending
an overconvergent isocrystal on a variety $X$ to a log-isocrystal
on a larger variety $\overline{X}$ is governed by the triviality of a sort of
local monodromy along components of the complement of $X$.
In \cite{kedlaya-part2}, it was shown that the problem can be localized on 
the space of valuations on the function field of the given variety.
In this paper, we solve the local semistable reduction problem at
monomial valuations (Abhyankar valuations of height 1 and residue transcendence
degree 0).

The context of this result (including a complex analogue) and a description
of potential applications is already given in the introduction of 
\cite{kedlaya-part1}, so we will not repeat it here. Instead, we devote
the remainder of this introduction to an overview of the results
specific to this paper, and a survey of the structure of the
various sections of the paper.

\subsection{Local semistable reduction}

The problem of (global) semistable reduction is to show that
an overconvergent $F$-isocrystal on a nonproper $k$-variety
can be extended to a log-isocrystal with nilpotent residues
on a proper $k$-variety
after passing to a suitable generically finite cover. In 
\cite{kedlaya-part1}, it was shown that existence of a log-extension
on a smooth pair $(X,D)$ can be checked generically along each
component of $D$. One can ensure the existence of the log-extension
along the proper transform of any given component, using the
$p$-adic local monodromy theorem of Andr\'e \cite{andre}, Mebkhout
\cite{mebkhout}, and this author \cite{kedlaya-local}; however,
when passing to the generically finite cover (and making sure it is
smooth by applying de Jong's alterations theorem \cite{dejong}), 
one typically introduces exceptional components along which one has not
achieved any control of local monodromy. See Example~\ref{exa:extra
monodromy} for an explicit
example of this phenomenon.

The main addition of \cite{kedlaya-part2} was to show that the  problem
of controlling exceptional components can be localized within
the Riemann-Zariski space of valuations of the function field of the
variety. Moreover, the resulting problem of local semistable reduction
can be reduced to a lower-dimensional case whenever one is working in
neighborhoods of a valuation which is composite (of height greater than 1),
or which has a residue field with
positive transcendence degree over the base field.

\subsection{Local monodromy at monomial valuations}

Monomial valuations on an $n$-dimensional variety
can be described as follows: in suitable local coordinates
$x_1, \dots, x_n$, they are determined by the fact that
$v(x_1), \dots, v(x_n)$ are linearly independent over $\QQ$.
On one hand, such valuations are particularly easy to describe,
so one expects to have an easier time working with them than with
other valuations. On the other hand, they form a subset of the Riemann-Zariski
space which, in some sense that we will not make precise here, is rather large.
(A related statement is that the set of Abhyankar places of a finitely
generated field extension is dense in the Riemann-Zariski space under the
patch topology \cite[Corollary~2]{kuhlmann}.)

In order to understand the structure of isocrystals in a neighborhood
of a monomial valuation, it is helpful to make some analysis
at the valuation itself. This
is the content of the paper \cite{kedlaya-fake}, which proves an analogue
of the $p$-adic local monodromy theorem for differential equations on a so-called
\emph{fake annulus} inside a higher-dimensional affine space.
It then makes sense to consider the (semisimplified) local monodromy 
representations attached to overconvergent $F$-isocrystals not just at
divisorial valuations, but also at monomial valuations. 

One is then led to ask how local monodromy varies as one varies
the monomial valuation, e.g., by varying $v(x_1), \dots, v(x_n)$.
We give a tangible answer to this question by defining
a higher-dimensional
analogue of the generic radius of convergence, as considered by Christol-Dwork
\cite{christol-dwork}. This gives a numerical invariant which
in the one-dimensional case computes the highest ramification break
of the local monodromy representation,
as in the work of Andr\'e, Christol-Mebkhout, Crew, Matsuda,
Tsuzuki, et al. (See \cite[\S~5]{kedlaya-mono-over} for an exposition.)
This number is shown to be a \emph{convex} function in $v(x_1), \dots, v(x_n)$
by the Hadamard three circles lemma in rigid geometry. 
(One can similarly construct an invariant that generalizes the full
Swan conductor in the one-dimensional case; 
see \cite{kedlaya-swan} for the beginning of this story.)

It is worth noting that this study has an interesting analogue over
the complex numbers, in the investigation of the
Stokes phenomenon conducted by Sabbah \cite{sabbah}. This concerns
irregular connections on complex surfaces, and (echoing an analogy
already seen in the one-dimensional situation) the variational behavior
of irregularity along divisors is apparently
quite similar to that of the invariant we consider.

Given what we have just described, we prove local semistable reduction
at a monomial valuation as follows. Using the $p$-adic local monodromy theorem
for fake annuli, we can force the highest ramification break
at the valuation itself to be zero. Then the convexity of the highest break
function implies that at certain nearby divisorial valuations, the highest
break is also forced to be zero, which forces the
isocrystal to be unipotent there also.

\subsection{Structure of the paper}

We conclude this introduction with a summary of the structure of the
paper.

In Section~\ref{sec:convex}, we derive some properties of convex functions.
The most important of these is a result which we were unable to find in 
the literature (Theorem~\ref{T:internally integral}), which shows that
a convex function whose values have the divisibility properties of a
piecewise affine function with integral coefficients must in fact be such
a function.

In Section~\ref{sec:diff alg}, we recall the relationship between
Newton polygons and norms of differential operators, as developed by
Christol-Dwork, Robba, Young, et al.

In Section~\ref{sec:gen radii}, we define our higher-dimensional analogue
of generic radius of convergence, and gather its key properties.

In Section~\ref{sec:around mono}, we introduce a form of the $p$-adic local
monodromy theorem covering so-called fake annuli. We then assert 
some related results,
notably the relationship between wild ramification and generic radius
of convergence for differential equations on $p$-adic curves.

In Section~\ref{sec:mono}, we develop some properties of monomial valuations,
then prove local semistable reduction at a
monomial valuation using the log-concavity of generic
radius of convergence.

In the Appendix, we describe two examples of semistable reduction.
One illustrates that one cannot insist on using a finite cover of a 
fixed compactification, rather than an alteration (as 
promised in the introduction of \cite{kedlaya-part1}). The other illustrates
that even if one starts with a good compactification, one cannot
achieve semistable reduction by doing so just for the divisors visible in
that compactification (as promised above).

\setcounter{theorem}{0}
\begin{notation}
We retain the basic notations of \cite{kedlaya-part1, kedlaya-part2}. In particular,
$k$ will always denote a field
of characteristic $p>0$, $K$ will denote a complete discretely valued
field of characteristic zero with residue field $k$, equipped with an
continuous 
endomorphism $\sigma_K$ lifting the $q$-power Frobenius for some power
$q$ of $p$, and $\gotho_K$ will denote the ring of integers of $K$.
\end{notation}

\subsection*{Acknowledgments}
Thanks to Yves Andr\'e for directing us to Sabbah's work.
Some of this material was presented at the Hodge Theory conference at
Venice International University in June 2006; that presentation 
was sponsored by the Clay Mathematics Institute.
The author was additionally supported by NSF grant DMS-0400727,
NSF CAREER grant DMS-0545904, and a Sloan Research Fellowship.

\section{Some properties of convex functions}
\label{sec:convex}

This section is completely elementary; it consists of some basic properties
of convex functions on subsets of $\RR^n$, which we will use later
to study variation of the highest ramification break as a function
of a valuation on a variety. We initially follow \cite{rock} for notation
and terminology.

\subsection{Convex functions}

In the study of convex functions, as in \cite{rock},
it is customary to use a slightly different setup than one might expect.
\begin{defn} \label{D:convex}
Denote $\RR_\infty = \RR \cup \{+\infty\}$.
A function $f: \RR^n \to \RR_\infty$
is \emph{convex} if for any $x,y \in \RR^n$ and $t \in [0,1]$, 
\[
tf(x) + (1-t)f(y) \geq f(tx + (1-t)y).
\]
Equivalently, $f$ is convex if and only if the \emph{epigraph} of $f$,
defined as 
\[
\epi(f) = \{(x_1, \dots, x_n, y) \in \RR^{n+1}: y \geq f(x_1,\dots,x_n)\},
\]
is a convex set.
\end{defn}

\begin{defn} \label{D:convex2}
If $U$ is a convex subset of $\RR^n$ and $f: U \to \RR$ is a function,
we say $f$ is \emph{convex} if the function $g: \RR^n \to \RR_\infty$ 
defined by
\[
g(x) = \begin{cases} f(x) & x \in U \\
+\infty & x \notin U
\end{cases}
\]
is convex in the sense of Definition~\ref{D:convex}.
Conversely, for $g: \RR^n \to \RR_\infty$ a convex function, we define the
\emph{essential domain} of $g$ to be 
\[
\dom(g) = \{x \in \RR^n: g(x) < +\infty\};
\]
then the restriction of $g$ to $\dom(g)$ is a convex function in the sense
just described. Write $\intdom(g)$ for the interior of $\dom(g)$;
then $g$ is continuous on $\intdom(g)$ \cite[Theorem~10.1]{rock}.
\end{defn}

\begin{defn}
For $C \subseteq \RR^n$, define the \emph{indicator function} $\delta_C: \RR^n
\to \RR_\infty$ by
\[
\delta_C(x) = \begin{cases} 0 & x \in C \\ +\infty & x \notin C;
\end{cases}
\]
then $\delta_C$ is a convex function if and only if $C$ is a convex set.
\end{defn}

\begin{defn}
An \emph{affine functional} is a map $\lambda: \RR^n \to \RR$ of the form
$\lambda(x_1,\dots,x_n) = a_1 x_1 + \cdots + a_n x_n + b$ for some
$a_1, \dots, a_n, b \in \RR$.
A \emph{generalized affine functional} is a map $\lambda: \RR^n \to \RR_\infty$
which is either an affine functional, or an affine functional plus
the indicator function of a \emph{closed halfspace}, i.e., a set of the 
form
\[
\{(x_1, \dots, x_n) \in \RR^n: a_1 x_1 + \cdots + a_n x_n \leq b\}
\]
with $a_1, \dots, a_n$ not all zero.
If $a_1, \dots, a_n, b \in \ZZ$ (in both places, if working in
the generalized case), 
we say
$\lambda$ is an \emph{integral} (generalized) affine functional.
\end{defn}

\begin{lemma} \label{L:inf convex}
Let $f: \RR^n \to \RR_\infty$ be a convex function. Then for any
$m \in \RR$, the function
$g_m: \RR^{n-1} \to \RR_\infty$ defined by
\[
g_m(x_1,\dots,x_{n-1}) = \inf_{x_n \in \RR} \{f(x_1, \dots, x_n) - mx_n\}
\]
is convex.
\end{lemma}
\begin{proof}
Since $f(x_1, \dots, x_n) - mx_n$ is again a convex function of
$x_1,\dots,x_n$, it suffices to consider the case $m=0$.
Given $x_{1,1},\dots,x_{1,n-1},x_{2,1},\dots,x_{2,n-1} \in \RR$
and $\epsilon > 0$,
choose $x_{1,n}, x_{2,n} \in \RR$ such that
\[
f(x_{i,1},\dots, x_{i,n})
\leq g_0(x_{i,1}, \dots, x_{i,n-1}) + \epsilon
\qquad (i=1,2).
\]
For $t \in [0,1]$, put $x_{3,j} = t x_{1,j} + (1-t)x_{2,j}$ for
$j=1,\dots,n$. 
Write $x_i = (x_{i,1},\dots,x_{i,n})$ and
$x'_i = (x_{i,1},\dots,x_{i,n-1})$ for $i=1,2,3$.
Then
\[
t g_0(x'_1) + (1-t) g_0(x'_2) + 2\epsilon
\geq
t f(x_1)
+ (1-t) f(x_2)
\geq f(x_3) \geq g_0(x'_3).
\]
Taking $\epsilon$ arbitrarily small, we deduce the convexity of $g_0$.
\end{proof}

\subsection{Internally polyhedral functions}

\begin{defn} \label{D:domain of aff}
For $f: \RR^n \to \RR_\infty$ a convex function, a \emph{domain of affinity}
is a subset $U$ of $\RR^n$ with nonempty interior on which $f$ agrees with
an affine functional $\lambda$. The nonempty interior condition
ensures that $\lambda$ is uniquely determined; we call it the
\emph{ambient functional} on $U$. 
\end{defn}

\begin{remark} \label{R:domain of aff}
Note that if $\lambda$ is an ambient functional on some
domain of affinity for $f$,
then the graph of $\lambda$ is a supporting hyperplane for the epigraph of $f$,
and so $f(x) \geq \lambda(x)$ for all $x \in \RR^n$.
\end{remark}

\begin{defn} \label{D:polyhedral}
A function $f: \RR^n \to \RR_\infty$
is \emph{polyhedral} if it has the form
\begin{equation} \label{eq:polyhedral}
f(x) = \max\{\lambda_1(x), \dots, \lambda_m(x)\}
\end{equation}
for some generalized affine functionals $\lambda_1,\dots,\lambda_m: \RR^n 
\to \RR_\infty$.
We say $f$ is \emph{integral polyhedral} if the $\lambda_i$ can be taken to
be integral; we say a set $C$ is \emph{rational polyhedral} if 
the indicator function $\delta_C$ is
integral polyhedral. 
We say $f$ is \emph{internally (integral) polyhedral} if
for each bounded (rational) polyhedral set $C \subseteq \intdom(f)$,
$f + \delta_C$ is (integral) polyhedral.
\end{defn}

\begin{remark}
It may look a bit strange to say that $C$ is \emph{rational} polyhedral if 
$\delta_C$ is \emph{integral} 
polyhedral. The point is that to get what one would properly call an
``integral polyhedral set'', i.e., the convex hull of a finite subset
of $\ZZ^n$, we would have to force the $\lambda_i$
in \eqref{eq:polyhedral} to have the form
$a_1 x_1 + \cdots + a_n x_n + b$ in which $a_1,\dots,a_n,b \in \ZZ$
but $a_1, \dots, a_n$ are additionally constrained to be coprime.
\end{remark}

\begin{remark}
The condition that a function $f$ be
 internally polyhedral is more permissive
than the condition that it be \emph{locally polyhedral}, in the sense of,
e.g., \cite[\S 15]{fujishige}. To say that $f$ is locally polyhedral means
that for every bounded polyhedral set $C$ meeting $\dom(f)$, $f + \delta_C$
is polyhedral. To see the difference, note that the functions
$f: \RR \to \RR_{\infty}$ given by
\[
f(x) = \begin{cases} +\infty & x \leq 0 \\
2N - N(N+1)x & x \in [1/(N+1), 1/N], N \in \ZZ_{>0} \\
0 & x \geq 1
\end{cases}
\]
and by
\[
f(x) = \begin{cases} +\infty & x < 0 \\
1 & x = 0 \\
0 & x > 0
\end{cases}
\]
are internally integral polyhedral but not locally polyhedral.
\end{remark}

\begin{lemma} \label{L:domains of aff}
\begin{enumerate}
\item[(a)] Let $f: \RR^n \to \RR_\infty$ be a convex function such that
$\dom(f)$ is polyhedral with nonempty interior. 
Then $f$ is polyhedral if and only if 
$\dom(f)$ is covered by finitely many domains of affinity for $f$.
\item[(b)] Let $f: \RR^n \to \RR_\infty$ be a convex function.
Then $f$ is internally polyhedral if and only if $\intdom(f)$ is covered by
(possibly infinitely many) domains of affinity for $f$.
\end{enumerate}
\end{lemma}
\begin{proof}
The ``only if'' implication is evident in both cases, so we focus on the
``if'' implications.

To prove (a), put $C = \dom(f)$; since $C$ is polyhedral,
we can write 
\[
\delta_C(x) = \max\{\mu_1(x), \dots, \mu_r(x)\}
\]
where each $\mu_j$ is the indicator function of a closed halfspace.
If $\dom(f)$ is covered by domains of affinity $U_1, \dots, U_m$
with ambient functionals $\lambda_1, \dots, \lambda_m$, then by
Remark~\ref{R:domain of aff}, we have
\[
f(x) = \max_{i,j}\{\lambda_i(x) + \mu_j(x)\},
\]
so $f$ is polyhedral. 

To prove (b), let $C$ be any bounded rational polyhedral subset of $\intdom(f)$.
Since $C$ is compact, it is covered by finitely many domains of affinity
for $f$; hence (a) implies that $f + \delta_C$ is polyhedral.
Since $\intdom(f)$ is the union of its bounded rational polyhedral subsets,
this proves the claim.
\end{proof}

\begin{cor} \label{C:infinite sup}
Let $f: \RR^n \to \RR_\infty$ be a convex function.
Then $f$ is internally integral polyhedral if and only if there exists
a (possibly infinite) subset $S$ of $\ZZ^n$ and a function $b: S \to \ZZ$
such that
\begin{equation} \label{eq:infinite sup}
f(x) = \sup_{s \in S} \{s_1 x_1 + \cdots + s_n x_n + b(s)\} \qquad
(x \in \intdom(f)).
\end{equation}
Moreover, in this case the supremum in \eqref{eq:infinite sup} is always
achieved.
\end{cor}
\begin{proof}
If $f$ is internally integral polyhedral, we choose $S$ and $b$
so that $s_1 x_1 + \cdots + s_n x_n + b(s)$ runs over the ambient
functionals on the domains of affinity of $f$; then
Remark~\ref{R:domain of aff} implies \eqref{eq:infinite sup}
with the supremum being achieved.
Conversely, suppose $S$ and $b$ exist.
Pick $x \in \intdom(f)$, and then choose $\epsilon > 0$ such that the box $B = 
\prod_{i=1}^n [x_i - 2\epsilon, x_i + 2\epsilon]$ is contained in
$\intdom(f)$.
Put $B' = \prod_{i=1}^n [x_i - \epsilon, x_i + \epsilon]$.
Let $U$ and $L$ be the supremum and infimum, respectively,
of $f$ on $B$ (which exist because $f$ is continuous on $\intdom(f)$).
If $s_1 > (U-L)/\epsilon$, then for $y = (y_1, \dots, y_n) \in B'$,
\begin{align*}
f(y) &\geq L \\
&> U - s_1 \epsilon \\
&\geq U - s_1(x_1 + 2\epsilon - y_1) \\
&= U - (s_1 (x_1 + 2\epsilon) + s_2 y_2 + \cdots + s_n y_n + b(s))
+ (s_1 y_1 + \cdots + s_n y_n + b(s)) \\
&\geq U - f(x_1+2\epsilon,y_2,\dots,y_n) + s_1 y_1 + \cdots + s_n y_n + b(s)\\
&\geq s_1 y_1 + \cdots + s_n y_n + b(s).
\end{align*}
That is, for all $x \in B'$, any term in \eqref{eq:infinite sup} for an $s$
with $s_1 > (U-L)/\epsilon$ can be omitted without changing the supremum.
Similarly, we can omit all $s$ for which $|s_i| > (U-L)/\epsilon$ for 
$i \in \{1,\dots,n\}$. Consequently,
we can compute the supremum in \eqref{eq:infinite sup} using only
finitely many affine functionals, and so $f + \delta_{C}$ is
integral polyhedral for any rational polyhedron $C \subseteq B'$.

Consequently, the point $x \in \intdom(f)$
admits a neighborhood contained in a union of finitely many
domains of affinity for $f$
corresponding to integral affine functionals. Since $x$ was arbitrary,
we deduce that all of $\intdom(f)$ can be covered by domains of affinity
for $f$ whose ambient functionals are integral.
By Lemma~\ref{L:domains of aff}, $f$ is internally integral polyhedral,
as desired.
\end{proof}

\subsection{Integral values and integral polyhedral functions}

The key result in this section (Theorem~\ref{T:internally integral})
asserts that the fact that a convex function is
internally integral polyhedral can be observed from its values at
rational $n$-tuples.

\begin{lemma} \label{L:base case}
Let $f: \RR^n \to \RR_\infty$ be a convex function such that
\begin{equation} \label{eq:intern int}
f(x_1, \dots, x_n) \in (\ZZ + \ZZ x_1 + \cdots + \ZZ x_n) \cup \{+\infty\}
\qquad (x_1, \dots, x_n \in \QQ).
\end{equation}
Then for any $x_1, \dots, x_{n-1} \in \QQ$, the function
$g: \RR \to \RR_\infty$ given by
$g(x) = f(x_1, \dots, x_{n-1}, x)$ is internally polyhedral, 
and on each domain of affinity of $g$, we have $g(x) = mx+ b$
for some $m \in \ZZ$ and some $b \in \ZZ + \ZZ x_1 + \cdots + \ZZ x_{n-1}$.
\end{lemma}
\begin{proof}
Fix $(x_1,\dots,x_{n-1}) \in \QQ^n$. We may assume that
$\intdom(g)$ is nonempty, as otherwise the claim is vacuously true; choose
$x_n \in \intdom(g)$.
Let $d$ be the least common
denominator of $x_1,\dots,x_n$, so that
\[
\frac{1}{d} \ZZ = \ZZ + \ZZ x_1 + \cdots + \ZZ x_n.
\]
For $N$ any sufficiently large positive integer, we have
\[
\frac{f(x_1,\dots,x_{n-1},x_n + 1/(dN)) - f(x_1,\dots,x_n)}{1/(dN)} 
\in dN \left(\ZZ + \ZZ x_1 + \cdots + \ZZ x_n + \frac{1}{dN} \ZZ \right)
= \ZZ.
\]
As $N \to \infty$,
this difference quotient runs through a sequence of integers which is nonincreasing
and bounded below (because $g$ is convex and $x_n \in \intdom(g)$).
Thus the quotient stabilizes for $N$ large. By convexity, the function
$g$ must be affine with
integral slope in a one-sided neighborhood of $x_n$;
since $\QQ$ is dense in $\RR$, the closed intervals on which $g$ is affine
with integral slope cover the interior of the essential domain. By
Lemma~\ref{L:domains of aff}, $g$ is internally polyhedral.

Let $d'$ be the least common denominator of $x_1, \dots, x_{n-1}$.
On any domain of affinity for $g$, we can write $g(x) = mx + b$ for some
$m \in \ZZ$. In this domain, we can find $y_1, y_2$ such that when we write
$y_i = r_i/s_i$ in lowest terms, we have $d',s_1,s_2$ coprime in pairs. 
From \eqref{eq:intern int}, we have $g(y_i) = m(r_i/s_i) + b \in 
\frac{1}{d' s_i} \ZZ$, implying $d' s_i b \in \ZZ$. Since this
holds for both $i=1$ and $i=2$, we find $d' b \in \ZZ$.
\end{proof}

\begin{theorem} \label{T:internally integral}
Let $f: \RR^n \to \RR_\infty$ be a convex function such that
\[
f(x_1, \dots, x_n) \in (\ZZ + \ZZ x_1 + \cdots + \ZZ x_n) \cup \{+\infty\}
\qquad (x_1, \dots, x_n \in \QQ).
\]
Then $f$ is internally integral polyhedral.
\end{theorem}
\begin{proof}
We proceed by induction on $n$, the case $n=1$ being solved by
Lemma~\ref{L:base case}. Write for brevity $x = (x_1, \dots, x_n)$ and
$x' = (x_1, \dots, x_{n-1})$.
For $m \in \ZZ$, define
\[
g_m(x') = \inf_{x \in \RR}
\{ f(x_1, \dots, x_{n-1}, x) - mx\};
\]
by Lemma~\ref{L:inf convex},
$g_m$ is a convex function on $\RR^{n-1}$.
By Lemma~\ref{L:base case},
for $x' \in \QQ^{n-1} \cap \intdom(g_m)$,
$g_m(x') \in \ZZ + \ZZ x_1 + \cdots + \ZZ x_{n-1}$.
We may thus apply the induction hypothesis to deduce that $g_m$
is internally integral polyhedral.

By one direction of Corollary~\ref{C:infinite sup}, we can construct
sets $S_m \subseteq \ZZ^{n-1}$ and functions $b_m: S_m \to \ZZ$ such that
\[
g_m(x') 
= \sup_{s \in S_m} \{s_1 x_1 + \cdots + s_{n-1} x_{n-1} + b_m(s)\}
\qquad (x' \in \intdom(g_m)).
\]
By Lemma~\ref{L:base case}, we know that
for $x \in \QQ^n$,
\[
f(x) = \sup_{m \in \ZZ} \{g_m(x_1,\dots,x_{n-1}) + mx_n\},
\]
so we conclude that
\begin{equation} \label{eq:final int rep}
f(x) = \sup\{s_1 x_1 + \cdots + s_{n-1} x_{n-1} + m x_n + b_m(s)\}
\end{equation}
for all $x \in \QQ^n \cap \intdom(f)$,
with the supremum running over $m \in \ZZ$ and $s = (s_1, \dots, s_{n-1}) \in S_m$.
Since both sides of \eqref{eq:final int rep} represent convex functions
on $\intdom(f)$ and they agree on a dense subset thereof,
we may invoke the continuity of convex functions to deduce that
\eqref{eq:final int rep} holds in fact for all $x \in
\intdom(f)$. 
By the other direction of Corollary~\ref{C:infinite sup}, 
$f$ is internally integral polyhedral, as desired.
\end{proof}

\subsection{Extension to rational polyhedral sets}

Although we will not use it in this paper, we note for future
reference a slight strengthening of Theorem~\ref{T:internally integral}.

\begin{lemma} \label{L:vertex finite}
Let $C$ be a bounded rational polyhedral subset of $\RR^n$, and let
$v \in \QQ^n$ be a vertex of $C$.
Let $f: \RR^n \to \RR_{\infty}$ be a convex function with $f(v) < \infty$,
and let $T_v$ be the set of integral affine functionals $\lambda$
which achieve their maximum on $C$ at $v$, and which agree with
$f$ on some domain of affinity meeting $C$. Then $T_v$ is finite.
\end{lemma}
\begin{proof}
Let $S \subset \ZZ^n$ be the set of $n$-tuples for which
\[
\max_{x \in C} \{ s_1 x_1 + \cdots + s_n x_n\}
= s_1 v_1 + \cdots + s_n v_n.
\]
Then $S$ is the intersection of $\ZZ^n$ with a strictly convex rational
polyhedral cone, and so is isomorphic to the intersection of $\ZZ^n_{\geq 0}$
with a sublattice of $\ZZ^n$ of finite index. Consequently, $S$
is well partially ordered, that is,
any infinite sequence of $S$
contains an infinite nondecreasing subsequence.

Suppose that $T_v$ is infinite. For $\lambda \in T_v$, write 
\[
\lambda(x) = s_1 x_1 + \cdots + s_n x_n + b
\]
with $s = s(\lambda) \in S$ and $b = b(\lambda)$.
Note that no two $\lambda \in T_v$ can have the same $s(\lambda)$,
by Remark~\ref{R:domain of aff}.
By the above,
we can choose $\lambda^{(1)}, \lambda^{(2)}, \dots \in T_v$ so that
the corresponding $s^{(i)} = s(\lambda^{(i)})$ form an infinite
increasing sequence.

For any $x$ in a domain of affinity for $f$ on which $f$ agrees with 
$\lambda^{(i)}$,
we must have $\lambda^{(i)}(x) = f(x)
\geq \lambda^{(i-1)}(x)$ because $f$ is convex;
moreover, because $\lambda^{(i)} \neq \lambda^{(i-1)}$,
we can choose $x$ so that the inequality is strict.
Since $s^{(i)} - s^{(i-1)} \in S$ by construction, we have
\[
\lambda^{(i)}(v) - \lambda^{(i-1)}(v)
\geq \lambda^{(i)}(x) - \lambda^{(i-1)}(x) > 0.
\]

However, for all $i$, we have $f(v) \geq \lambda^{(i)}(v)$. Consequently,
the $\lambda^{(i)}(v)$ form a strictly increasing, bounded above sequence
with values in the discrete subset $\ZZ + \ZZ v_1 + \cdots + \ZZ v_n$
of $\RR$. This is impossible, yielding the desired contradiction.
\end{proof}

\begin{theorem} \label{T:integral2}
Let $C$ be a bounded rational polyhedral
subset of $\RR^n$. Then a continuous
convex function $f: C \to \RR$
is integral polyhedral if and only if 
\begin{equation} \label{eq:integral values}
f(x) \in \ZZ + \ZZ x_1 + \cdots + \ZZ x_n \qquad (x \in C \cap \QQ^n).
\end{equation}
\end{theorem}
\begin{proof}
If $f$ is integral polyhedral, then \eqref{eq:integral values} is clear.
Conversely, assume \eqref{eq:integral values}; by
Theorem~\ref{T:internally integral}, $f$ is internally integral polyhedral.

Let $T$ be the set of integral affine functionals $\lambda$ which agree
with $f$ on some domain of affinity.
For $v$ a vertex of $C$, let $T_v$ be the set of $\lambda \in T$ which
achieve their maximum on $C$ at $v$.
By Lemma~\ref{L:vertex finite}, each $T_v$ is finite;
since $C$ has only finitely many vertices, and $T$ is the union of the
$T_v$, $T$ must be finite.
By Lemma~\ref{L:domains of aff}, $f$ is integral polyhedral, as desired.
\end{proof}

\section{Differential equations and Newton polygons}
\label{sec:diff alg}

In this section, we review the relationship between differential equations
over complete valued fields and Newton polygons. The analysis here 
draws from Young \cite{young-tams},
Christol-Dwork \cite{christol-dwork}, and particularly
Robba \cite{robba}.

\setcounter{theorem}{0}
\begin{hypothesis} \label{H:valued diff}
Until \S\ref{subsec:higher},
let $F$ denote a \emph{valued (nontrivial) differential field}
of characteristic zero. That is, $F$ 
is a field equipped with a nonzero derivation $\del: F \to F$, and with
a nonarchimedean absolute value $|\cdot|$; we write
$v(\cdot) = - \log |\cdot|$ for the corresponding valuation.
We will later require that $F$ be complete
(starting in \S\ref{subsec:splitting}).
\end{hypothesis}

\subsection{Valued differential fields and twisted polynomials}

\begin{defn}
For $T$ a bounded linear operator on a normed vector space $V$, the \emph{operator
norm} of $T$, denoted $|T|_V$, is the infimum of 
those $c \in \RR_{\geq 0}$ for which
$|T(x)| \leq c|x|$ for all $x \in V$. For $m,n \in \ZZ_{\geq 0}$,
we have the evident inequality
\[
|T^{m+n}|_V \leq |T^m|_V |T^n|_V.
\]
By taking logarithms, we arrive at
the situation of Fekete's lemma: if $\{a_n\}_{n=1}^\infty$ is a sequence of
reals with
$a_{m+n} \leq a_m + a_n$ for all $m,n$, then 
the sequence $\{ a_n/n\}_{n=1}^\infty$ 
either converges to its infimum or diverges to $-\infty$
\cite[Part 1, Problem~98]{polya-szego}.
We may thus define the \emph{spectral norm} of $T$ as
\[
|T|_{V,\spect} = \lim_{n \to \infty} |T^n|_V^{1/n} = \inf_n \{|T^n|_V^{1/n}\};
\]
it depends only on the equivalence class of the norm on $V$.
In particular, we will apply this notation with $T = \del$ acting on $F$
(as a vector space over the subfield killed by $\del$);
put $r_0 = -\log |\del|_F$.
\end{defn}

\begin{defn}
Let $F\{T\}$ denote the twisted polynomial ring over $F$ in the sense of Ore
\cite{ore}, so that for $x \in F$, $Tx = xT + \del(x)$.
By the Leibniz rule, for $x \in F$, 
\[
T^n x = \sum_{i=0}^n \binom{n}{i} \del^{n-i}(x) T^i.
\]
The twisted polynomial ring admits division with remainder on both sides,
so the Euclidean algorithm applies to show that left ideals and right
ideals are all principal (again, see \cite{ore}).
\end{defn}

\begin{remark} \label{R:opposite}
Note that the opposite ring of $F\{T\}$ is also a twisted polynomial ring,
for the difference field $F'$ with the same underlying field as $F$, but
with derivation $-\del$. The passage to the opposite ring corresponds, in
the classical language of 
differential equations, to replacing a differential
operator with its adjoint.
\end{remark}

\begin{defn}
For $P = \sum c_n T^n \in F\{T\}$, define the \emph{Newton polygon} of $P$
as the lower convex hull of the set
\[
\{(-n, v(c_n)): n \in \ZZ_{\geq 0}, c_n \neq 0\}.
\]
Define the \emph{multiplicity} of a real number $r$ 
(as a slope of $P$)
as the 
width of the segment of the Newton polygon of slope $r$, or 0 if there is
no such segment.
For $r \in \RR$, define
\[
v_r(P) = \min_n \{rn + v(c_n)\};
\]
this is the $y$-intercept of the supporting line of the Newton polygon
of slope $r$. Note that for $P$ fixed, $v_r(P)$ 
is a continuous function of $r$.
\end{defn}

As originally observed by Robba \cite[\S 1]{robba},
this Newton polygon behaves like its counterpart for untwisted polynomials, but
only for slopes which are not too large.
\begin{lemma} \label{L:add diff slopes}
For $P,Q \in F\{T\}$ and $r \leq r_0$, we have
$v_r(PQ) = v_r(P) + v_r(Q)$.
\end{lemma}
\begin{proof}
Write $P = \sum_i a_i T^i$ and $Q = \sum_j b_j T^j$; then 
\[
PQ = \sum_k \left( \sum_{i+j=k} \sum_{h\geq 0}
\binom{i+h}{h} a_{i+h} \del^h(b_j) \right) T^k,
\]
and hence
\begin{equation} \label{eq:Newton slope}
\begin{split}
v_r(PQ) &\geq \min_{h,i,j} \{ v(a_{i+h}) + v(b_j) + r(i+j) - \log |\del^{h}|_F\} \\
&\geq \min_{h,i,j} \{ v(a_{i+h}) + v(b_j) + r(i+j) + h r_0 \} \\
&\geq \min_{h,i,j} \{ v(a_{i+h}) + v(b_j) + r(i+h+j)\}.
\end{split}
\end{equation}
This immediately yields $v_r(PQ) \geq v_r(P) + v_r(Q)$. To  establish 
equality for $r < r_0$, 
let $i_0$ and $j_0$ be the smallest values of $i$ and $j$
which minimize 
$ri + v(a_i)$ and $rj + v(b_j)$, respectively; then \eqref{eq:Newton slope}
achieves its minimum for $h=0, i = i_0, j=j_0$ but not for any other
$h,i,j$ with $i+j = i_0 + j_0$. Hence $v_r(PQ) = v_r(P) + v_r(Q)$; equality
for $r = r_0$ follows by continuity.
(Compare \cite[Proposition~1.6(2)]{robba}.)
\end{proof}
\begin{cor} \label{C:slopes add}
For $P,Q \in F\{T\}$ and $r < r_0$, the multiplicity of $r$
as a slope of $PQ$ is the sum of its multiplicities as a slope of $P$
and of $Q$.
\end{cor}

The moral here is that when one is only looking at
phenomena in slopes less than $r_0$, one does not see the
difference between twisted and untwisted polynomials. For instance,
here is an explicit instance of this conclusion modeled on
\cite[Lemme~1.4]{christol-dwork} (compare also
\cite[Proposition~1.6(1)]{robba}).

\begin{lemma} \label{L:untwist norm}
For $r \leq r_0$,
let $Q = U^d + \sum_{i=0}^{d-1} a_i U^i \in F[U]$
be a polynomial with all slopes at least $r$.
(Here $F[U]$ denotes the untwisted polynomial ring.)
Put $W = F[U]/F[U]Q$ as an $F$-vector space with norm
$|\sum_{i=0}^{d-1} c_i U^i| = \max\{|c_i| e^{-ri}\}$.
Let $U$ act on $W$ by left multiplication, and let $\del$ act
coordinatewise with respect to the basis $1,U,\dots,U^{d-1}$.
Then 
\[
|(U+\del)^n|_W \leq e^{-rn} \qquad \mbox{for all $n \in \ZZ_{\geq 0}$};
\]
moreover, equality holds in case
$r < r_0$
and $Q$ has all slopes equal to $r$.
\end{lemma}
\begin{proof}
Rewrite the slope hypothesis as
$|a_i|_F \leq e^{-r(d-i)}$
for $i=0, \dots, d-1$;
then clearly $|U^n|_W \leq |U|_W^n = e^{-rn}$, so
$|(U + \del)^n - U^n|_W \leq e^{-r(n-1)} |\del|_F$. This
yields all of the claims.
\end{proof}

\subsection{Splitting over a complete field}
\label{subsec:splitting}

For $F$ complete, we obtain Robba's analogue for differential operators
\cite[Th\'eor\`eme~2.4]{robba}
of Hensel's lemma for an untwisted polynomial over a complete
nonarchimedean field.
\begin{hypothesis} \label{H:complete}
Throughout this subsection and the next, 
assume that $F$ is complete for its norm.
\end{hypothesis}

\begin{prop} \label{P:diff factor}
Fix $r < r_0$ and $m \in \ZZ_{\geq 0}$.
Let $R \in F\{T\}$ be a twisted polynomial such that 
$v_r(R - T^m) > v_r(T^m)$. Then $R$ can be factored uniquely as $PQ$, where
$P \in F\{T\}$ has degree $\deg(R) - m$ and
all slopes less than
$r$, $Q \in F\{T\}$ is monic of degree $m$ and has all slopes
greater than $r$, $v_r(P - 1) > 0$, and $v_r(Q - T^m) > v_r(T^m)$.
\end{prop}
\begin{proof}
We first check existence. Define sequences $\{P_l\}, \{Q_l\}$ as follows.
Define $P_0 = 1$ and $Q_0 = T^m$. Given $P_l$ and $Q_l$, write
\[
R - P_l Q_l = \sum_i a_i T^i,
\]
then put
\[
X_l = \sum_{i \geq m} a_i T^{i-m}, \qquad
Y_l = \sum_{i < m} a_i T^i
\]
and set $P_{l+1} = P_l + X_l$, $Q_{l+1} = Q_l + Y_l$.
Put $c_l = v_r(R- P_l Q_l) - rm$, so that $c_0 > 0$.
Suppose that $v_r(P_l-1) \geq c_0$, $v_r(Q_l-T^m) \geq c_0 + rm$,
and $c_l \geq c_0$.
Then visibly $v_r(P_{l+1}-1) \geq c_0$ and $v_r(Q_{l+1}-T^m) \geq c_0 + rm$;
by Lemma~\ref{L:add diff slopes},
\begin{align*}
c_{l+1} &= v_r(R - (P_l + X_l)(Q_l + Y_l)) - rm \\
&= v_r(X_l(T^m - Q_l) + (1-P_l)Y_l - X_lY_l) - rm \\
&\geq \min\{c_l + (c_0 +rm), c_0 + (c_l+rm) , c_l + (c_l+rm)\} - rm \\
&\geq c_l + c_0.
\end{align*}
By induction on $l$, we deduce that $c_l \geq (l+1) c_0$.
Moreover, each $P_l$ has degree at most $\deg(R)-m$,
and each $Q_l - T^m$ has degree at most $m-1$.
Consequently, the sequences $\{P_l\}$ and $\{Q_l\}$ converge 
under $v_r$ to polynomials $P$ and $Q$, which 
have the desired properties.

We next check uniqueness. Suppose $R = P_1 Q_1$ is a second such factorization;
put $c = \min\{v_r(P - P_1), v_r(Q - Q_1) - v_r(T^m)\}$.
Put
\[
X = R - P_1 Q = (P - P_1)Q = P_1(Q_1 - Q),
\]
and suppose $X \neq 0$; then $v_r(X) = c + rm$
by Lemma~\ref{L:add diff slopes}.
Write $X = \sum b_k T^k$, and choose $k$ such that
$v_r(X) = v_r(b_k T^k)$. The equality
\[
X = (P - P_1)T^m + (P - P_1)(Q - T^m)
\]
shows that we cannot have $k < m$, while the equality
\[
X = Q_1 - Q + (P_1 - 1)(Q_1 - Q)
\]
shows that we cannot have $k \geq m$. This contradiction
forces $X=0$, proving $P = P_1, Q = Q_1$ as desired.
\end{proof}

\begin{remark}
Note that the proof of Proposition~\ref{P:diff factor} does not
involve any divisions. Consequently,
if the coefficients of $P$ lie in a subring $S$ of $F$
which is complete under the norm, then the coefficients of $Q$ and $R$
also lie in $S$.
\end{remark}

We obtain a corollary akin to a factorization result of Dwork-Robba
\cite[Theorem~6.2.3]{dwork-robba}.
\begin{cor} \label{C:slope factor}
Any monic twisted polynomial $P \in F\{T\}$ admits a unique factorization
\[
P = P_{r_1} \cdots P_{r_m} P_+
\]
for some $r_1 < \cdots < r_m < r_0$,
where each $P_{r_i}$ is monic with all slopes equal to $r_i$,
and $P_+$ is monic with all slopes at least $r_0$.
\end{cor}

\begin{remark} \label{R:opposite factor}
Note that by Remark~\ref{R:opposite}, Corollary~\ref{C:slope factor} can
also be stated with the factors in the reverse order;
the degrees of the individual factors will not change, but the
factors themselves may differ.
\end{remark}

\subsection{Differential modules}

Remember that we are still assuming that $F$ is complete 
(Hypothesis~\ref{H:complete}).

\begin{defn} \label{D:diff module}
A \emph{differential module} over $F$ is a finite dimensional
$F$-vector space $V$ equipped with an action of $\del$, or equivalently,
a left $F\{T\}$-module which is finite over $F$.
Given a basis $B$ of $V$, we may equip $V$ with the supremum norm with 
respect to $B$, and thus define operator and spectral norms $|\del|_{V,B}$
and $|\del|_{V,B,\spect}$. Changing $B$ gives an equivalent norm on $V$, so the
spectral norm $|\del|_{V,B,\spect}$ does not depend on $B$; we thus write it
as $|\del|_{V,\spect}$.
\end{defn}

\begin{remark}
We will also have occasion to speak about differential modules over
fields equipped with multiple derivations, in which case the notation
for the operator/spectral norm will indicate which derivation is being
measured. See \S\ref{subsec:higher}.
\end{remark}

\begin{remark} \label{R:diff module}
Instead of the spectral norm $|\del|_{V,\spect}$, we will
invariably consider the truncated spectral norm
$\max\{|\del|_{F,\spect}, |\del|_{V,\spect}\}$.
(It turns out that these coincide \cite[Lemma~6.2.4]{kedlaya-course},
but we will not use that fact here.)
The truncated spectral norm
can be computed in terms of a basis of $V$ as follows:
if $D_n$ denotes the matrix
via which $\del^n$ acts on this basis, then
\begin{equation} \label{eq:diff module basis}
\max\{|\del|_{F,\spect}, |\del|_{V,\spect}\}
= \max\{|\del|_{F,\spect}, \limsup_{n \to \infty} |D_n|^{1/n} \},
\end{equation}
where the norm applied to $D_n$ is the supremum over entries
\cite[Proposition~1.3]{christol-dwork}.
\end{remark}

\begin{defn}
Let $V$ be a differential module over $F$. 
A \emph{cyclic vector} for $V$ is an element $\bv \in V$ not contained in
any proper differential submodule; it is equivalent to ask
that $\bv, \del(\bv),
\dots, \del^{n-1}(\bv)$ form a basis of $V$ for $n = \dim_F(V)$.
A cyclic vector defines an isomorphism $V \cong F\{T\}/F\{T\}P$ for
some $P \in F\{T\}$. 
\end{defn}

\begin{lemma} \label{L:cyclic vector}
Every differential module over $F$ contains a 
cyclic vector.
\end{lemma}
\begin{proof}
See, e.g., \cite[Theorem~III.4.2]{g-functions}.
\end{proof}

\begin{lemma}  \label{L:exact seq}
Let $P \in F\{T\}$ be a monic twisted polynomial and let
$V = F\{T\}/F\{T\}P$ be the corresponding differential module. Then
every short exact sequence $0 \to V_1 \to V \to V_2 \to 0$ of
differential modules arises uniquely from a factorization $P = P_1 P_2$
of $P$ into monic twisted polynomials, in which $V_1 \cong F\{T\}/F\{T\}P_1$
and $V_2 \cong F\{T\}/F\{T\}P_2$ and the map $V \to V_2$ is induced by the
natural projection $F\{T\}/F\{T\}P \to F\{T\}/F\{T\}P_2$.
\end{lemma}
\begin{proof}
The kernel of $F\{T\} \to V_2$ is a left ideal of $F\{T\}$, so it is generated
by a unique monic $P_2$, giving the isomorphism $V_2 \cong F\{T\}/F\{T\}P_2$
and the factorization $P = P_1P_2$. We also have
$V_1 \cong F\{T\}P_2 / F\{T\}P$, and the latter is isomorphic to
$F\{T\}/F\{T\}P_1$ via right division by $P_2$.
\end{proof}

The following is attributed ``Dwork-Katz-Turritin'' (sic) in 
\cite[Th\'eor\`eme~1.5]{christol-dwork}.
\begin{theorem} \label{T:slope to norm}
Let $P \in F\{T\}$ be a nonzero twisted polynomial with least slope $r$,
and put $V = F\{T\}/F\{T\}P$. Then
\[
\max\{|\del|_F, |\del|_{V,\spect}\} = \max\{|\del|_F, e^{-r}\}.
\]
\end{theorem}
\begin{proof}
If $P$ has a single slope $r$ and that slope satisfies $r < r_0$, or 
if $P$ has all slopes at least $r_0$, then we obtain the claim by using
the basis $1, T, \dots, T^{\deg(P)-1}$ and invoking
Lemma~\ref{L:untwist norm}.
Otherwise, we may apply Corollary~\ref{C:slope factor} to reduce to such cases.
\end{proof}
\begin{remark}
The proof of \cite[Th\'eor\`eme~1.5]{christol-dwork} contains a minor
error 
in its implication 1$\implies$2: in its notation, one passes from
$K$ to an algebraic extension $K(z)$ without worrying about whether $\|D\|$
increases as a result. (In our notation, this amounts to passing
from $F$ to an extension
without checking whether $|\del|_F$ increases.)
The proof of Theorem~\ref{T:slope
to norm} shows that the final result is nonetheless correct,
and indeed the proof is only slightly changed.
\end{remark}

\begin{lemma} \label{L:one slope}
Let $P \in F\{T\}$ be a nonzero twisted polynomial with all slopes
equal to $r < r_0$ (resp.\ all slopes at least $r_0$). 
Then every Jordan-H\"older factor $W$
of $V = F\{T\}/F\{T\}P$ satisfies $|\del|_{W,\spect} = e^{-r}$
(resp. $|\del|_{W,\spect} \leq |\del|_F$).
\end{lemma}
\begin{proof}
We induct on $\dim_F(V)$. If $V$ is irreducible, then
Theorem~\ref{T:slope to norm} implies the claim. Otherwise,
choose a short exact sequence $0 \to V_1 \to V \to V_2 \to 0$;
by Lemma~\ref{L:exact seq},
we have a factorization $P = P_1 P_2$ such that
$V_i \cong F\{T\}/F\{T\}P_i$ for $i=1,2$.
By Corollary~\ref{C:slopes add},
$P_1$ and $P_2$ also have all slopes equal to $r$
(resp.\ all slopes at least $r_0$),
so we may apply the induction hypothesis to $V_1, V_2$ to conclude.
\end{proof}

\begin{theorem} \label{T:norm decomp}
Let $V$ be a differential module over $F$. 
Then there is a unique decomposition
\[
V = V_+ \oplus \bigoplus_{r < r_0} V_r
\]
of differential modules, such that each Jordan-H\"older factor $W_+$
of $V_+$ satisfies $|\del|_{W_+,\spect} \leq |\del|_F$, and
each Jordan-H\"older factor $W_r$ of $V_r$ satisfies
$|\del|_{W_r,\spect} = e^{-r}$.
\end{theorem}
\begin{proof}
The decomposition is clearly unique if it exists. To produce it, 
we induct on $\dim_F(V)$.
Choose a cyclic vector, let $V \cong F\{T\}/F\{T\}P$ be the resulting
isomorphism, and let $r_1$ be the least slope of $P$. If
$r_1 \geq r_0$, we may take $V = V_+$ and be done, so assume 
$r_1 < r_0$. If $P$ has all slopes equal to $r_1$, then
Lemma~\ref{L:one slope} implies that we may take $V = V_{r_1}$
and be done, so assume the contrary.

Apply Corollary~\ref{C:slope factor} to factor
$P = P_{r_1} Q$ with $P_{r_1}$ having all slopes equal to $r_1$,
and $Q$ having all slopes greater than $r_1$.
By Lemma~\ref{L:exact seq}, this factorization gives rise to an exact sequence
$0 \to V_1 \to V \to V_2 \to 0$ in which (by Lemma~\ref{L:one slope}
and the induction hypothesis) each Jordan-H\"older factor of $V_1$
has spectral norm of $\del$ equal to 
$e^{-r_1}$, and each Jordan-H\"older factor of $V_2$
has spectral norm of $\del$ strictly less than $e^{-r_1}$.

Now apply Corollary~\ref{C:slope factor} to factor $P$ again, but this time
in the opposite ring of $F\{T\}$ as per Remark~\ref{R:opposite factor}. 
That is, write $P = Q' P'_{r_1}$ with $P'_{r_1}$ having all slopes equal to $r_1$
and $Q'$ having all slopes greater than $r_1$.
Then Lemma~\ref{L:exact seq} and Lemma~\ref{L:one slope} give
an exact sequence
$0 \to V'_1 \to V \to V'_2 \to 0$ in which each Jordan-H\"older factor
of $V'_2$ has spectral norm $e^{-r_1}$, and each Jordan-H\"older factor of
$V'_1$ has spectral norm strictly less than $e^{-r_1}$.
In particular, $\dim(V_1) = \dim(V'_2)$, and
$V_1 \cap V'_1 = \{0\}$; this forces $V \cong V_1 \oplus V'_1$.
Splitting off
$V_{r_1} = V_1$ and repeating, we obtain 
the desired decomposition.
\end{proof}

\subsection{Differential fields of higher order}
\label{subsec:higher}

\begin{hypothesis} \label{H:valued diff2}
We now modify Hypothesis~\ref{H:valued diff}
to say that $F$ is a \emph{complete valued differential field of order $n$}
of characteristic zero. That is, in addition to being complete for
a norm, $F$ is equipped
with not one but $n$ commuting nonzero
derivations $\del_1, \dots, \del_n$.
\end{hypothesis}

When comparing norms for different derivations acting on a differential module, 
it is useful to renormalize
to remove the spectral norms of the derivations themselves.
\begin{defn}
A \emph{differential module} over $F$ is now a finite dimensional
$K$-vector space $V$ equipped with actions of $\del_1, \dots, \del_n$.
Define the \emph{scale} of $V$ as
\[
\max\left\{\max\left\{1, \frac{|\del_i|_{V,\spect}}{|\del_i|_{F, \spect}}\right\}: i = 1, \dots, n \right\}.
\]
For each $i$ at which the outer maximum is achieved, we say $\del_i$ is \emph{dominant for $V$}.
\end{defn}

\begin{theorem} \label{T:scale decomp}
Suppose that the $|\del_i|_F/|\del_i|_{F,\spect}$ for $i=1, \dots, n$
are all equal to a common value $s_0$. 
Let $V$ be a differential module over $F$. Then
there is a unique decomposition
\[
V = V_+ \oplus \bigoplus_{s > s_0} V_s
\]
such that each Jordan-H\"older factor of $V_s$ has scale $s$,
and each Jordan-H\"older factor of $V_+$ has scale at most $s_0$.
\end{theorem}
\begin{proof}
Apply Theorem~\ref{T:norm decomp} for each $\del_i$; the uniqueness 
assertion in the proposition means that the decomposition
with respect to $\del_i$ is respected by the
other $\del_j$. By taking the common refinement of these decompositions,
then appropriately recombining terms, we obtain the desired result.
\end{proof}

\begin{prop} \label{P:diff rational}
Suppose that $F$ is discretely valued, and 
the $|\del_i|_F/|\del_i|_{F,\spect}$ for $i=1, \dots, n$
are all equal to a common value $s_0$. 
Let $V$ be a differential module over $F$, and consider the
decomposition in Theorem~\ref{T:scale decomp}. 
Then for $s > s_0$, $s^{\dim(V_s)} \in s_0^{\ZZ} |F^*|$.
\end{prop}
\begin{proof}
Apply Theorem~\ref{T:norm decomp} to $V_s$ for each $\del_i$. From
the result, we obtain a decomposition $V_s = \oplus V_{s,i}$ in which
for each Jordan-H\"older factor $W_i$ of $V_{s,i}$,
we have that $\del_j$ is dominant for $W_i$ 
when $j=i$ but not when $j =1,\dots,i-1$.

Choose a cyclic vector for $V_{s,i}$ with respect to $\del_i$; 
let $P(T) = T^d + \sum_{i=0}^{d-1} a_i T^i$ be the resulting twisted polynomial.
By Corollary~\ref{C:slope factor}, the Newton polygon of $P$ must
have all slopes equal to $-\log (s |\del_i|_{F,\spect})$; it follows that
$(s |\del_i|_{F,\spect})^d = |a_0|$.

Note that $|\del_i|_F \in |F^*|$ because $F$ is discretely valued, so $|\del_i|_{F,\spect} = |\del_i|_F/s_0
\in s_0^\ZZ |F^*|$, and so $s^{\dim(V_{s,i})} \in s_0^\ZZ |F^*|$.
Since $\sum_i \dim(V_{s,i}) = \dim(V_s)$,
$s^{\dim(V_s)} \in s_0^\ZZ |F^*|$, as desired.
\end{proof}

\section{Generic radii of convergence}
\label{sec:gen radii}

In this section, we revisit the usual notion of generic radii of convergence
of differential equations from the work of Dwork, Robba, et al., but this time
working in several dimensions.

\subsection{Generalized polyannuli}

It will be convenient to consider subsets of affine spaces more general
than the polyannuli considered in \cite[Definition~3.1.5]{kedlaya-part1}.

\begin{notation}
For $X = (X_1, \dots, X_n)$ an $n$-tuple:
\begin{itemize}
\item
for $A$ an $n \times n$ matrix, write $X^A$ for the $n$-tuple
whose $j$-th entry is $\prod_{i=1}^n x_i^{A_{ij}}$;
\item
for $B$ an $n$-tuple, put $X^B = X^A$ for $A$ the diagonal matrix with $A_{ii} = B_i$;
\item
for $c$ a number, put $X^c = X^A$ for $A$ the scalar matrix $cI_n$.
\end{itemize}
\end{notation}

\begin{defn}
By a \emph{log-(rational polyhedral) subset},
or \emph{log-RP subset}, of $(0, +\infty)^n$, we will 
mean a subset
$S$ whose image under the logarithm map to $\RR^n$ is a rational
polyhedral set in the sense of Definition~\ref{D:polyhedral}.
We say $S$ is \emph{ind-log-RP} if it is the union of an increasing
sequence of log-RP subsets.
\end{defn}

\begin{notation}
Let $S$ be an ind-log-RP subset of $(0, +\infty)^n$.
Write $A_K(S)$ for the rigid analytic subspace of $\AAA^n_K$ defined by the conditions
\[
(|t_1|, \dots, |t_n|) \in S;
\]
if $S$ is log-RP and $\log(S)$ is bounded, then $A_K(S)$ is affinoid.
Note that
\[
\Gamma(A_K(S), \calO) = 
\left\{
\sum_{J \in \ZZ^n} c_J T^J: c_J \in K,
\lim_{J \to \infty} |c_J| R^J = 0
\quad (R \in S) \right\},
\]
where $T = (t_1, \dots, t_n)$. (The limit condition should be
interpreted as follows: for each $R \in S$
and each $\epsilon>0$, there are only finitely many $J \in \ZZ^n$ with
$|c_J| R^J > \epsilon$.)
For $S = \{R\}$ a singleton set, we write $A_K(R)$ for $A_K(S)$.
\end{notation}

The following toric coordinate changes will be useful.
\begin{defn}
For $A$ an $n \times n$ matrix,
let $f_A: \GG_m^n \to \GG_m^n$ be the map $T \mapsto T^A$,
or the induced map $A_K(R) \to A_K(R^A)$.
\end{defn}
\begin{lemma} \label{L:toric}
For any complete extension $K'$ of $K$, 
any
$c_1, \dots, c_n \in K'$ with $|c_i| = r_i$,
and any $\lambda \in (0,1]$,
define the open polydisc
\[
D(C,\lambda R) = \{T \in A_K(R): |t_i - c_i| < \lambda r_i
\quad (i=1, \dots, n)\}.
\]
Then $f_A$ carries $D(C,\lambda R)$ to $D(C^A, \lambda R^A)$.
\end{lemma}
\begin{proof}
Rewrite the defining condition of $D(C, \lambda R)$ as
$|1 - t_i/c_i| < \lambda$. Then note that this implies
$|1 - (t_i/c_i)^n| < \lambda$ for any $n \in \ZZ$,
by examination of the binomial expansion of $(1 - (1 - t_i/c_i))^n$.
To finish, recall that for $\lambda \in (0,1]$, $|1-a|, |1-b| < \lambda$
implies $|1-ab| < \lambda$ because $1-ab = (1-a) + (1-b) - (1-a)(1-b)$.
\end{proof}

\begin{defn}
For $R = (r_1, \dots, r_n) \in S$, the space $A_K(S)$
carries a Gauss norm $|\cdot|_R$ defined by
\[
\left|\sum_J c_J T^J \right|_R = \sup_J \{|c_J| R^J \};
\]
it is in fact the supremum norm on $A_K(R)$.
\end{defn}

The following convexity lemma (analogous to the Hadamard three circles theorem)
is a repackaging of \cite[Lemma~3.1.6]{kedlaya-part1}, but similar observations
occur much earlier in the literature, e.g., 
\cite[Corollaire~4.2.8]{amice}, \cite[Corollaire~5.4.9]{christol-robba}.
\begin{lemma} \label{L:three circles}
For $A,B \in S$ and $c \in [0,1]$, put $R = A^c B^{1-c}$; that is,
$r_i = a_i^c b_i^{1-c}$ for $i=1, \dots, n$.
Then for any $f \in \Gamma(A_K(S), \calO)$,
\[
|f|_R \leq |f|_A^c |f|_B^{1-c}.
\]
\end{lemma}
\begin{proof}
Since each Gauss norm
is calculated as a supremum over monomials, it suffices
to check the inequality in the case of a single monomial, 
in which case it becomes an equality. 
\end{proof}

\subsection{Generic radii of convergence}

\begin{defn} \label{D:radius}
Let $S$ be a log-RP subset of $(0,+\infty)^n$, take $R \in S$, and
let $\calE$ be a $\nabla$-module (locally free coherent sheaf plus
integrable connection) on $A_K(S)$. 
Let $F$ (or $F_R$ in case of ambiguity)
be the completion of $\Frac \Gamma(A_K(S),\calO)$ under $|\cdot|_R$, and 
put 
\[
V = \Gamma(A_K(S), \calE) \otimes_{\Gamma(A_K(S),\calO)}
 F.
\]
For $i=1, \dots, n$, define $\del_i = \frac{\del}{\del t_i}$ as a
derivation on $F$. View $F$ as a differential field of order $n$, 
view $V$ as a differential module over $F$, and let $T(\calE,R)$
be the reciprocal of the scale of $V$; that is,
\[
T(\calE, R) = \min_i \{ \min\{1, |\del_i|_{V,\spect}^{-1} |\del_i|_{F,\spect}\}\}.
\]
\end{defn}

\begin{remark} \label{R:geometric radius}
We may interpret $T(\calE,R)$ as the largest $\lambda \in (0,1]$ such that
for any complete extension $K'$ of $K$ and any $C = (c_1, \dots, c_n) \in (K')^n$
with $|c_i| = r_i$ for $i=1, \dots, n$, $\calE$ admits a basis of
horizontal sections on $D(C,\lambda R)$.
In particular, for $n=1$, our function $T(\calE, R)$ is equal to $R^{-1}$ times the
generic radius of convergence $R(\calE,R)$ of \cite{christol-dwork}.
The letter $T$ is used here to denote ``toric normalization''.
\end{remark}

\begin{remark}
It may be helpful to compare Remark~\ref{R:geometric radius}
with \cite[Definition~5.3]{kedlaya-mono-over},
but one must beware of three typos in the latter:
the $\min$ should be a $\max$, the subscript $\rho$ is missing, 
and the reference to \cite[Proposition~1.2]{christol-dwork} should be
to Proposition 1.3 therein.
\end{remark}

\begin{remark} \label{R:radius ops}
The following are easily verified.
\begin{itemize}
\item
If $0 \to \calE_1 \to \calE \to \calE_2 \to 0$ is exact, then
\[
T(\calE, R) = \min\{T(\calE_1,R), T(\calE_2,R)\}.
\]
\item
We have
\[
T(\calE_1 \otimes \calE_2, R)\leq \min\{T(\calE_1,R), T(\calE_2,R)\}.
\]
\item
We have
\[
T(\calE^\dual, R) = T(\calE,R).
\]
\end{itemize}
\end{remark}

The function $T$ also satisfies a toric invariance property.
\begin{prop} \label{P:toric}
Let $S$ be a log-RP subset of $(0,+\infty)^n$, take $R \in S$, and
let $\calE$ be a $\nabla$-module on $A_K(S)$. 
For $A \in M_n(\ZZ)$, put $S^A = \{R^A: R \in S\}$.
Then for any $\nabla$-module $\calE$ on $A_K(S^A)$,
\[
T(f_A^* \calE, R) \geq T(\calE,R^A),
\]
with equality if $A \in \GL_n(\ZZ)$.
\end{prop}
\begin{proof}
This follows immediately from Lemma~\ref{L:toric}.
\end{proof}

Lemma~\ref{L:three circles} yields the 
following log-concavity property, which  generalizes
\cite[Proposition~2.3]{christol-dwork}.
\begin{prop} \label{P:log concave}
Let $S$ be a log-RP subset of $(0,+\infty)^n$.
Let $\calE$ be a $\nabla$-module on $A_K(S)$.
For $A,B \in S$ and $c \in [0,1]$,
\[
T(\calE,A^c B^{1-c}) \geq T(\calE,A)^c T(\calE,B)^{1-c}.
\]
\end{prop}
\begin{proof}
Since $S$ is log-RP, $A_K(S)$ is affinoid; by Kiehl's theorem,
$\calE$
is generated by finitely many global sections. Let $\be_1,\dots,\be_m$
be a maximal linearly independent set of global sections, and let
$D_{i,l}$ be the matrix over $\Frac \Gamma(A_K(S), \calO)$ via which
$\frac{\del^l}{\del t_i^l}$ acts on $\be_1, \dots, \be_m$. Since
$\be_1, \dots, \be_m$ are maximal linearly independent, we can choose
$f \in \Gamma(A_K(S), \calO)$ so that $f \Gamma(A_K(S), \calE)$
is contained in the span of $\be_1, \dots, \be_m$. This implies that
$f D_{i,l}$ has entries in $\Gamma(A_K(S), \calO)$ for all $i,l$.

Put $R = A^c B^{1-c}$.
By Lemma~\ref{L:three circles}, we have
\[
|f|_R |D_{i,l}|_R \leq |f|^c_A |D_{i,l}|^c_A |f|^{1-c}_B |D_{i,l}|^{1-c}_B;
\]
taking $l$-th roots of both sides and taking limits superior
yields
\begin{align*}
\max\{|\del_i|_{F_R, \spect},
\limsup_{l \to \infty} |D_{i,l}|_R^{1/l}\} &\leq
\max\left\{ |\del_i|_{F_A, \spect},
\left( \limsup_{l \to \infty} |D_{i,l}|_A^{1/l} \right)^c \right\} \\
&\qquad \cdot
\max\left\{ |\del_i|_{F_B, \spect},
\left( \limsup_{l \to \infty} |D_{i,l}|_B^{1/l} \right)^{1-c} \right\}
\end{align*}
because the factors coming from $f$ all tend to 1.
By \eqref{eq:diff module basis}, this yields the desired result.
\end{proof}

\begin{example} \label{exa:dwork}
Let $\calE$ be the $\nabla$-module of rank 1 defined by
$\nabla \bv = \lambda \pi d(t_1^{i_1} \cdots t_n^{i_n})$, where 
$\lambda \in \gotho_K^*$,
$\pi \in K$ satisfies $\pi^{p-1} = -p$ (that is, $\pi$ is a
\emph{Dwork pi} and $\calE$ is a \emph{Dwork isocrystal}),
and $i_1,\dots,i_n \in \ZZ$ are not all divisible by $p$. Then
as in \cite[Chapter 5]{kedlaya-mono-over}, one may check
that
\[
T(\calE,R) = \min\{1, r_1^{-i_1} \cdots r_n^{-i_n}\}.
\]
\end{example}

\subsection{The Frobenius antecedent theorem}

We now revisit the Frobenius antecedent theorem of Christol-Dwork
\cite[Th\'eor\`eme~5.4]{christol-dwork} in a higher-dimensional context,
following \cite[Theorem~6.15]{kedlaya-mono-over}.

\begin{hypothesis} \label{H:frob ant}
Let $Y$ be an affinoid space over $K$, and suppose $t_1, \dots, t_n \in 
\Gamma(Y, \calO)^*$
are such that $dt_1, \dots, dt_n$ freely generate $\Omega^1_{Y/K}$;
let $f: Y \to \AAA^n_K$ be the resulting \'etale morphism.
Form the Cartesian diagram
\begin{equation} \label{eq:cartesian}
\xymatrix{
Y' \ar^g[r] \ar[d] & Y \ar^f[d] \\
\AAA^n_K \ar[r] & \AAA^n_K
}
\end{equation}
in which the morphism $\AAA^n_K \to \AAA^n_K$ is given by
$t_i \mapsto t_i^p$ $(i=1,\dots,n)$.
Let $\calE'$ be a $\nabla$-module on $Y'$ such that
\begin{equation} \label{eq:Frob condition}
\left| \frac{\del}{\del t_i} \right|_{\calE', \spect}
< |t_i|^{-1}_{\sup, Y} \qquad (i=1,\dots,n),
\end{equation}
where the left-hand side is computed using any norm on 
$\Gamma(Y',\calE')$ compatible with the affinoid norm on
$\Gamma(Y',\calO)$. (Since any two such norms are equivalent, the
spectral norm does not depend on the choice.)
\end{hypothesis}

\begin{defn} \label{D:action}
Suppose that $K$ contains a primitive $p$-th root of unity $\zeta$.
For $J = (j_1, \dots, j_n) \in (\ZZ/p\ZZ)^n$,
let $g_J: Y' \to Y'$ be the map defined by $t_i \mapsto t_i \zeta^{j_i}$ 
for $i=1, \dots, n$. (More precisely, we get $g_J$ from the 
Cartesian square \eqref{eq:cartesian} using the original map
$Y' \to Y$ and the map $Y' \to \AAA^n_K$ given by $t_1 \zeta^{j_1},
\dots, t_n \zeta^{j_n}$.)
Then the map $h_J: g_J^* \calE' \to \calE'$ defined by
\[
h_J(\bv) = \sum_{i_1,\dots,i_n=0}^\infty (\zeta^{j_1} - 1)^{i_1}
\cdots (\zeta^{j_n}-1)^{i_n} \frac{t_1^{i_1} \cdots t_n^{i_n}}{i_1!\cdots i_n!}
\frac{\del^{i_1}}{\del t_1^{i_1}} \cdots
\frac{\del^{i_n}}{\del t_n^{i_n}} \bv
\]
converges 
because of \eqref{eq:Frob condition}.
\end{defn}

\begin{prop} \label{P:frob ant1}
Suppose that $K$ contains a primitive $p$-th root of unity $\zeta$.
Under Hypothesis~\ref{H:frob ant},
there is a unique $\nabla$-module $\calE$ on $Y$ such that
$g^* \calE \cong \calE'$ and
the action of the $h_J$ on $\calE'$ is induced by the trivial
action on $\calE$.
\end{prop}
\begin{proof}
Put $M' = \Gamma(Y',\calE')$.
The maps $h_J$ satisfy $h_J(t_i \bv) = \zeta^{j_i} t_i h_J(\bv)$; hence
for $J = (j_1, \dots, j_n) \in \{0, \dots, p-1\}^n$, if we define
\[
f_J(\bv) = t_1^{-j_1} \cdots t_n^{-j_n} 
\sum_{J' \in (\ZZ/p\ZZ)^n}
\zeta^{-j_1j'_1-\cdots-j_n j'_n} h_{J'}(\bv),
\]
then $f_J(\bv)$
is fixed by the $h_{J'}$. 
Let $M$ be the $\Gamma(Y', g^{-1}(\calO))$-span of the
$f_J(\bv)$; then $M$ is a coherent $\Gamma(Y,\calO)$-module,
and (by an appropriate form of Hilbert's Theorem 90)
the natural map $M \otimes \Gamma(Y',\calO) \to M'$ is a
$(\ZZ/p\ZZ)^n$-equivariant isomorphism. We give $M$ a
$\nabla$-module structure
by declaring the action of $\frac{\del}{\del t_i}$ on $M$
to be $p^{-1} t_i^{1-p}$ times the action of
$\frac{\del}{\del t_i}$ on $M'$.
This gives rise to $\calE$ such that
$\calE' \cong g^* \calE$, which evidently is unique for the property
of being fixed by the $h_J$.
\end{proof}

\begin{defn}
Under Hypothesis~\ref{H:frob ant}, we call $\calE'$ the
\emph{Frobenius antecedent} of $\calE$. Note that the uniqueness
implies that it makes sense to define a Frobenius antecedent 
for a $\nabla$-module on a rigid space $Y$ even if 
\eqref{eq:Frob condition}
is only satisfied after replacing $Y$ with each element of an admissible
open cover,
 or if $K$ does not contain a primitive $p$-th root of unity.
\end{defn}

\subsection{Frobenius antecedents and generic radii}

\begin{notation}
Throughout this subsection, 
write $S^{1/p} = \{R^{1/p}: R \in S\}$ for $S \subseteq
(0, +\infty)^n$, and let $f_p$ denote the map $f_{pI_n}: A_K(S^{1/p}) \to A_K(S)$ 
for any $S$.
\end{notation}

\begin{lemma} \label{L:toric2}
Let $S$ be a log-RP subset of $(0,+\infty)^n$, suppose
$(1,\dots,1,\rho) \in S$, and
let $\calE$ be a $\nabla$-module on $A_K(S)$. 
Then
\[
T(f_{p}^* \calE,(1,\dots,1,\rho^{1/p}))
\geq 
T(\calE, (1,\dots,1,\rho))^{1/p}.
\]
\end{lemma}
This inequality can be shown to be an equality when
$T(\calE,(1,\dots,1,\rho)) > |p|^{p/(p-1)}$, but we will not need that 
more precise result here.

\begin{proof}
It suffices to observe that 
for $R = (1,\dots,1,\rho)$,
$f_{p}$ carries $D(C^{1/p}, \lambda^{1/p}
R^{1/p})$ into $D(C, \lambda R)$.
The latter follows from \cite[Lemma~5.12]{kedlaya-mono-over},
but note a misprint therein: in the last line of the statement,
the quantities $r\rho^{1/p}$ and $r^p \rho$ should be
$r^{1/p} \rho^{1/p}$ and $r\rho$, respectively.
\end{proof}

\begin{theorem} \label{T:frob ant}
Put $S = [1,1]^{n-1} \times (\epsilon,1)$ for some $\epsilon \in (0,1)$.
Let $\calF$ be a $\nabla$-module on $A_K(S^{1/p})$ such that
\begin{equation} \label{eq:frob ant}
T(\calF,(1, \dots, 1, \rho^{1/p})) > |p|^{1/(p-1)} \qquad
(\rho \in (\epsilon, 1)).
\end{equation}
Then $\calF$ admits a Frobenius antecedent $\calE$ on $A_K(S)$,
which satisfies
\begin{equation} \label{eq:convergence0}
T(\calE, (1,\dots,1,\rho))^{1/p} = T(\calF,(1,\dots,1,\rho^{1/p})) 
\qquad (\rho \in (\epsilon,1)).
\end{equation}
\end{theorem}
\begin{proof}
For each point in $S$, \eqref{eq:frob ant}
implies that one can find a neighborhood $S'$ of that point in $S$ such that
on $A_K(S')$, \eqref{eq:Frob condition} holds.
We then glue to obtain a Frobenius antecedent on all of $A_K(S)$.

To prove \eqref{eq:convergence0}, note that
with $R = (1,\dots,1,\rho)$,
given $c_1, \dots, c_n$ with $|c_i| = 1$ for $i=1,\dots,n-1$
and $|c_n| = \rho^{1/p} \in (\epsilon^{1/p},1)$, we can apply the maps $f_J$
(from the proof of Proposition~\ref{P:frob ant1}) to horizontal
sections on a polydisc $D(C, \lambda^{1/p} 
R^{1/p})$ to obtain horizontal sections
on $D(C^p, \lambda R)$. Consequently,
\[
T(\calE, (1,\dots,1,\rho))^{1/p} \geq T(\calF,(1,\dots,1,\rho^{1/p})) 
\qquad (\rho \in (\epsilon,1));
\]
the reverse inequality follows from Lemma~\ref{L:toric2}.
\end{proof}

Using Frobenius antecedents, one overcomes the scale barrier
built into the results of Section~\ref{sec:diff alg}.
\begin{lemma} \label{L:integral}
Take $S$ as in Theorem~\ref{T:frob ant},
and let $\calE$ be a $\nabla$-module on $A_K(S)$.
Then for each $\rho \in (\epsilon,1)$, there exists
an integer $j \in \{1, \dots, \rank(\calE)\}$ 
and a nonnegative integer $m$ such that
\[
T(\calE,(1,\dots,1,\rho))^j \in \rho^\ZZ (|K^*| |p|^{(1/(p-1)) \ZZ})^{p^{-m}}.
\]
\end{lemma}
\begin{proof}
Let $m$ be the least nonnegative integer such that
\[
T(\calE, (1, \dots, 1, \rho))^{p^m} \leq |p|^{1/(p-1)}.
\]
If $T(\calE, (1, \dots, 1, \rho))^{p^m} = |p|^{1/(p-1)}$, then
we are done, so assume not. By Proposition~\ref{P:log concave},
$T(\calE,(1,\dots,1,\rho))$ is a log-concave and hence continuous function
of $\rho$, so we can choose a closed interval $I$ with endpoints in the
divisible closure of $|K^*|$, such that 
$T(\calE,(1,\dots,1,\eta)) > |p|^{p^{1-m}/(p-1)}$ for $\eta \in I$.
Apply Theorem~\ref{T:frob ant} $m$ times to produce a
$\nabla$-module $\calE'$ with
\[
T(\calE,(1,\dots,1,\rho)) = T(\calE',(1,\dots,1,\rho^{p^m}))^{1/p^m}.
\]
Then apply Proposition~\ref{P:diff rational} to $\calE'$,
noting that for the derivation
$\frac{\del}{\del t_i}$ on $\Frac \Gamma(A_K(S), \calO)$
under the $R$-Gauss norm, the operator norm and spectral norm are
$r_i^{-1}$ and $|p|^{1/(p-1)} r_i^{-1}$, respectively.
This yields the desired result.
(Compare \cite[Th\'eor\`eme~4.2-1]{cm3}.)
\end{proof}

\begin{lemma} \label{L:rational}
Take $S$ as in Theorem~\ref{T:frob ant},
and let $\calE$ be a $\nabla$-module on $A_K(S)$.
Suppose that $T(\calE, (1, \dots, 1, \rho)) \to 1$ as $\rho \to 1^-$.
Then there exist $\eta \in [\epsilon,1)$,
an integer $1 \leq j \leq \rank(\calE)$, and
a nonnegative integer $i$ such that
$T(\calE, (1, \dots, 1, \rho)) = \rho^{i/j}$ for $\rho \in (\eta,1)$.
\end{lemma}
\begin{proof}
There is no harm in assuming that $|p|^{1/(p-1)} \in |K^*|$.
For $c \in (0, -\log(\epsilon))$, define
\[
f(c) = \log T(\calE, (1, \dots, 1, e^{-c}));
\]
this function is concave by Proposition~\ref{P:log concave},
takes nonpositive values, and by hypothesis has
limit 0 as $c \to 0^+$. Consequently, $f$ is nonincreasing.

For $i$ a sufficiently large integer, we can find $c_i \in (0, -\log(\epsilon))$
such that 
\[
f(c_i) = \frac{1}{p^m(p-1)} \log |p|;
\]
the $c_i$ then form a decreasing sequence.
By Lemma~\ref{L:integral}, for each 
$c \in (c_{i+1}, c_i) \cap \QQ \log |p|$,
there exists $j \in \{1, \dots, \rank(\calE)\}$ such that
\[
f(c) \in \frac{1}{j}(p^{-m-1} \log |K^*| + \ZZ c).
\]
By Theorem~\ref{T:internally integral}, $f$ is piecewise affine on 
$(c_{i+1},c_i)$,
and each slope is a rational number with denominator
bounded by $\rank(\calE)$. In particular, the slopes of $f$ belong to a discrete
subgroup of $\RR$. 

As $c \to 0^+$, the slopes of $f$ on successive domains of affinity
form a nondecreasing sequence of values, each of which is nonpositive
because $f$ is nonincreasing. Since these values lie in a discrete subgroup
of $\RR$, they must stabilize; that is, $f$ is affine in some neighborhood of 0.
Since $f \to 0$ as $c \to 0^+$, $f$ must actually be linear in a neighborhood
of 0. This yields the desired result.
(Compare  \cite[Th\'eor\`eme~4.2-1]{cm3}.)
\end{proof}

\begin{defn}
We say an $n$-tuple $R \in (0,+\infty)^n$ is \emph{commensurable}
if $r_1, \dots, r_n$ generate a discrete subgroup of the multiplicative
group $\RR_{>0}$.
In this case, we call the generator of that subgroup lying in
$(0,1)$ the \emph{generator} of $R$.
\end{defn}

\begin{theorem} \label{T:rational}
Let $R  \neq (1,\dots,1) \in (0, +\infty)^n$ be commensurable with generator $\rho$.
Let $S$ be an ind-log-RP subset of $(0, +\infty)^n$ 
containing $R^c$ for all $c >0$ sufficiently
small. Let $\calE$ be a $\nabla$-module on $A_K(S)$ such that
$T(\calE,R^c) \to 1$ as $c \to 0^+$. Then there exist integers $i,j$ with
$i \geq 0$ and $1 \leq j \leq \rank(\calE)$, such that
\[
T(\calE, R^c) = \rho^{i/j} \qquad \mbox{for $c>0$ sufficiently small}.
\]
\end{theorem}
\begin{proof}
This reduces to Lemma~\ref{L:rational} by applying a suitable toric change of
coordinates $f_A$.
\end{proof}

\begin{remark} \label{R:Frobenius}
As in the one-dimensional case
\cite[Proposition~6.3-11]{cm3}, one can enforce the
condition that $T(\calE,R^c) \to 1$ as $c \to 0^+$ by equipping $\calE$
with a Frobenius structure. Explicitly, suppose that $q$ is a power of $p$,
and that 
$\sigma: A_K(S^{1/q}) \to A_K(S)$ is a map obtained by composing the toric map
$f_{qI_n}$ with a $q$-power Frobenius lift on $K$.
If there is an isomorphism $\sigma^* \calE \cong \calE$ over $A_K(S^{1/q})$,
then Lemma~\ref{L:toric2} implies
that for $R \in S$, $T(\calE,R^{1/{q^m}}) \geq T(\calE, R)^{1/q^m}$,
so the values of $T(\calE,R^c)$ get arbitrarily close to 1;
by Proposition~\ref{P:log concave}, it follows that
$T(\calE,R^c) \to 1$ as $c \to 0^+$.
\end{remark}

\section{Around the local monodromy theorem}
\label{sec:around mono}

In this section, we recall the $p$-adic local monodromy theorem, in a
generalized form suited to treating monomial valuations.
We then mention some related results, on the interplay between generic radii
of convergence in the one-dimensional case and local monodromy.

\subsection{The monodromy theorem for fake annuli}

To state the monodromy theorem at the level of generality we
need, we must recall some terminology from \cite{kedlaya-fake}.

\begin{defn}
We say a linear functional $\lambda: \RR^n \to \RR$
is \emph{irrational} if $\ZZ^n \cap \ker(\lambda) = \{0\}$. For $\lambda$
an irrational functional, write
$\lambda_1, \dots, \lambda_n$ for the images under $\lambda$ of
the standard generators of $\ZZ^n$.
For $I \subseteq (0,1)$, let $\calR^\lambda_I$ (resp.\
$\calR^{\lambda,\inte}_I$)
be the Fr\'echet completion of $K[t_1^{\pm}, \dots, t_n^{\pm}]$
(resp. $\gotho_K[t_1^{\pm}, \dots, t_n^{\pm}]$)
with respect to the Gauss norms $|\cdot|_{\rho^{\lambda_1}, \dots, \rho^{\lambda_n}}$
for $\rho \in I$. Write $\calR^\lambda$ (resp.\ $\calR^{\lambda,\inte}$)
for the union of
$\calR^\lambda_{[\rho,1)}$ 
(resp.\ $\calR^{\lambda,\inte}_{[\rho,1)}$)
over all $\rho \in (0,1)$. 
\end{defn}

\begin{remark}
In our notation for generalized polyannuli, $\calR^\lambda_I$ would be the
global sections of the structure sheaf on $A_K(S)$ for 
\[
S 
= \{(\rho^{\lambda_1}, \dots, \rho^{\lambda_n}): \rho \in I\}
\]
if the latter were an ind-log-RP subset; 
however, that can only happen when
$I$ consists of a single point, or when $n=1$ (the case of a true annulus).
This is what is fake about a so-called fake annulus: it does not 
fit any conventional definition of an analytic subspace of $\AAA^n_K$,
even in Berkovich's framework for nonarchimedean analytic
geometry \cite{berkovich}.
\end{remark}

\begin{remark}
Given an interval $I$, 
let $I'$ be the interval consisting of those $r \in (0, +\infty)$ such that
$|p|^{r/w(p)} \in I$. 
For $\lambda$ an irrational functional,
the ring $\Gamma^\lambda_{I'}$ of 
\cite[Definition~2.4.1]{kedlaya-fake} (with the lattice therein taken to be
$\ZZ^n$) is isomorphic to $\calR^\lambda_{I}$ via a map
sending $\{z_i\}$ to $t_i$ for $i=1, \dots, n$.
This identification has a number of consequences, some captured in
Lemma~\ref{L:consequences} below.
\end{remark}
\begin{lemma} \label{L:consequences}
\begin{enumerate}
\item[(a)]
For $I$ closed, $\calR^\lambda_{I}$ is a principal ideal
domain.
\item[(b)]
For any $\rho \in (0,1)$,
$\calR^\lambda_{[\rho,1)}$ is a B\'ezout domain (an integral domain
whose finitely generated ideals are principal).
\item[(c)]
Let $I_1 \subset I_2 \subset \cdots$ be an increasing sequence of closed intervals
with union $[\rho,1)$. Given any sequence $M_1, M_2, \dots$ in which $M_l$ is a
finite free $\calR^\lambda_{I_l}$-module, together with isomorphisms
$\iota_l: M_{l+1} \otimes \calR^\lambda_{I_l} \cong M_l$, there exist
a finite free $\calR^\lambda_{[\rho,1)}$-module $M$ and isomorphisms
$\psi_l: M \otimes \calR^\lambda_{I_l} \cong M_l$ such that
$\iota_l \circ \psi_{l+1} = \psi_l$; moreover, $M$ and the $\psi_l$ are
determined up to unique isomorphism.
\end{enumerate}
\end{lemma}
\begin{proof}
For (a), see \cite[Proposition~2.6.8]{kedlaya-slope}.
For (b), see \cite[Theorem~2.9.6]{kedlaya-slope}.
For (c), see \cite[Theorem~2.8.4]{kedlaya-slope}.
\end{proof}

\begin{defn}
Define a \emph{$\nabla$-module} over $\calR^\lambda$ as a
finite free $\calR^\lambda$-module $M$ equipped with an integrable connection
$\nabla: M \to M \otimes \Omega^1_{\calR^\lambda/K}$.
We say a $\nabla$-module over $\calR^\lambda$ is \emph{constant} if it
has a basis of horizontal sections, \emph{quasi-constant} if it
becomes constant after tensoring with a finite \'etale extension of
$\calR^{\lambda,\inte}$, and \emph{(quasi)-unipotent} if it admits a
filtration by $\nabla$-submodules whose successive quotients are
(quasi)-constant.
\end{defn}

\begin{defn}
Let $\sigma: \calR^\lambda \to \calR^\lambda$ be a 
continuous endomorphism lifting a power of the absolute
Frobenius map on the residue field
of $\calR^{\lambda, \inte}$.
Define an \emph{$F$-module} (resp. $(F,\nabla)$-module)
over $\calR^\lambda$ relative to $\sigma$
as a finite free $\calR^\lambda$-module (resp.\ $\nabla$-module)
$M$ equipped with an isomorphism
$F: \sigma^* M \to M$ of modules (resp.\ of $\nabla$-modules).
As with true annuli,
the category of $(F, \nabla)$-modules over $\calR^\lambda$
is canonically independent
of the choice of $\sigma$ \cite[Proposition~3.4.7]{kedlaya-fake}.
\end{defn}

\begin{defn}
For $s = c/d \in \QQ$, an $F$-module $M$ is \emph{pure} (or \emph{isoclinic})
\emph{of slope $s$}
if there exists a basis of $M$ on which $F^d$ acts via the product of
a scalar of valuation $c$ with an invertible matrix over
$\calR^{\lambda, \inte}$. Note that this is the equivalent characterization
of \cite[Proposition~6.3.5]{kedlaya-slope} rather than the original definition;
one can in fact develop the slope theory for $F$-modules using this definition
instead, as in \cite{kedlaya-rel}.
\end{defn}

In this language, one has the following result from \cite{kedlaya-fake}.

\begin{theorem} \label{T:fake annuli local mono}
Let $\calE$ be an $(F, \nabla)$-module over $\calR^\lambda$.
\begin{enumerate}
\item[(a)]
There exists a unique filtration $0 \subset \calE_1 \subset
\cdots \subset \calE_m \subset \calE$ of $\calE$ by 
$(F, \nabla)$-submodules such that each $\calE_i/\calE_{i-1}$ is
pure of some slope $s_i$ as an $F$-module, 
and $s_1 < \cdots < s_m$.
\item[(b)]
Each successive quotient of the filtration in (a) is quasi-constant
as a $\nabla$-module.
Consequently, $\calE$ is quasi-unipotent as a $\nabla$-module.
\end{enumerate}
\end{theorem}
\begin{proof}
Statement (a) is \cite[Theorem~5.2.1]{kedlaya-fake}; note that this depends on the
generalized slope filtration theorem of \cite{kedlaya-slope}, not just on the
original form of the theorem of \cite{kedlaya-local}. Statement (b) is
\cite[Theorem~5.2.4]{kedlaya-fake}.
\end{proof}

\subsection{Monodromy and convergence (one-dimensional case)}

We now revert from fake annuli back to true annuli, to recall some
results relating generic radii of convergence to wild ramification.
We defer to \cite{kedlaya-mono-over} for a more extensive discussion of the points
we only summarize here, including attributions.

\begin{notation}
Throughout this subsection, we take $n=1$, drop $\lambda$, and write $t$ for $t_1$.
Also, 
as in \cite{kedlaya-part1}, when we write an interval $I$ out explicitly, we typically
omit the parentheses in the notation $A_K(I)$.
\end{notation}

\begin{prop} \label{P:rep category}
The category of quasi-unipotent $\nabla$-modules over $\calR$ is equivalent to
the category of representations of
\[
\Gal(k((t))^{\sep}/k((t))) \times K
\]
in finite dimensional $K^{\unr}$-vector spaces, 
which are semilinear and permissible (the restriction to some open subgroup
is trivial) on the first factor, and algebraic, $K$-rational, and unipotent on the second factor.
\end{prop}
\begin{proof}
See \cite[Theorem~4.45]{kedlaya-mono-over}.
\end{proof}

\begin{defn}
Let $\calE$ be a $\nabla$-module on $A_K(\epsilon,1)$ for some $\epsilon \in
(0,1)$. By Lemma~\ref{L:consequences}(c), $\calE$ corresponds to a
$\nabla$-module over $\calR_{(\epsilon,1)}$; let
$M_\calE$ be the corresponding $\nabla$-module over $\calR$.
\end{defn}
\begin{prop} \label{P:ramif break}
Assume that the field $k$ is perfect.
Let $\calE$ be a $\nabla$-module on $A_K(\epsilon,1)$ for some $\epsilon \in
(0,1)$ such that $M_\calE$ is quasi-unipotent.
Then for $\rho \in (0,1)$ sufficiently close to $1$,
$T(\calE,\rho) = \rho^\beta$ for $\beta$ equal to the highest ramification break
of the Galois factor of the representation associated to $M_\calE$
by Proposition~\ref{P:rep category}. Moreover, if $\beta>0$ and
the lowest ramification
break is also equal to $\beta$, then for $\rho \in (0,1)$ sufficiently close 
to $1$, \emph{every} nonzero 
local horizontal section of $\calE$ around a generic
point of radius $\rho$ has exact radius of convergence $\rho^{\beta+1}$.
\end{prop}
\begin{proof}
See \cite[Theorem~5.23]{kedlaya-mono-over}.
\end{proof}

\begin{cor} \label{C:close break}
Let $\calE$ be a $\nabla$-module on $A_K(\epsilon,1)$ for some $\epsilon \in
(0,1)$, such that $M_\calE$ is quasi-unipotent. Then the following are equivalent.
\begin{enumerate}
\item[(a)]
There exists a positive integer $m$ coprime to $p$ such that
$M_\calE \otimes \calR[t^{1/m}]$ is unipotent.
\item[(b)]
$T(\calE,\rho) = 1$ for $\rho \in (\epsilon,1)$ sufficiently close to $1$.
\item[(c)]
$T(\calE,\rho) > \rho^{1/\rank(\calE)}$ for $\rho \in (\epsilon,1)$ sufficiently close to $1$.
\end{enumerate}
\end{cor}
\begin{proof}
There is no harm in enlarging $K$, so we may assume $k$ is perfect.
Clearly (a)$\implies$(b)$\implies$(c). Given (c), 
by Proposition~\ref{P:ramif break},
the highest ramification
break of the corresponding Galois representation must be less than 
$1/\rank(\calE)$; since the highest break must be a 
nonnegative rational number
with denominator at most $\rank(\calE)$ (by the Hasse-Arf theorem), it
must equal 0, that is, the 
representation is only tamely ramified. This
yields the claim.
\end{proof}

\subsection{Monodromy and convergence (relative case)}

In light of Proposition~\ref{P:ramif break}, it is natural to make the following definition.
\begin{defn}
With notation as in Theorem~\ref{T:rational}, we call the rational number $i/j$ the
\emph{(differential) highest ramification break} of $\calE$ in the direction of $R$,
denoted $b(\calE,R)$.
\end{defn}

\begin{prop} \label{P:log concave break}
Let $A,B \in (0, +\infty)^n$ be commensurable, take $c \in [0,1] \cap \QQ$,
put $R = A^c B^{1-c}$, and suppose $R$ is also commensurable.
Let $\alpha,\beta,\rho$ be the generators of $A,B,R$, respectively.
Let $S$ be a ind-log-RP subset of $(0, +\infty)^n$ which
contains $A^h, B^h, R^h$ for 
$h >0$ sufficiently small. Let $\calE$ be a $\nabla$-module on $A_K(S)$ such that
$T(\calE,*^h) \to 1$ as $h \to 0^+$ for $* \in \{A,B,R\}$. 
Then
\[
\rho^{b(\calE,R)} \geq \alpha^{c b(\calE,A)} \beta^{(1-c) b(\calE,B)}.
\]
\end{prop}
\begin{proof}
Apply Proposition~\ref{P:log concave}.
\end{proof}

\begin{defn} \label{D:generic fibre}
Take $S$ as in Theorem~\ref{T:frob ant}, and let $\calE$ be a $\nabla$-module on $A_K(S)$.
Let $L$ be the completion of $K(t_1, \dots, t_{n-1})$ under the $(1,\dots,1)$-Gauss norm.
Let $\calE$ be a $\nabla$-module on $A_K(S)$.
Let $I_1 \subset I_2 \subset \cdots$ be an increasing sequence of
closed intervals with union $(\epsilon,1)$.
Put $S_l = [1,1]^{n-1} \times I_l$, and put
\[
M_l = \Gamma(A_K(S_l), \calE) \otimes_{\Gamma(A_K(S_l), \calO)} 
\Gamma(A_L(I_l), \calO);
\]
then there is a unique locally free coherent sheaf $\calF$ on
$A_L(\epsilon,1)$ admitting identifications $M_l \cong \Gamma(A_L(I_l), \calF)$
compatible with restriction. Moreover, $\calF$ inherits the structure
of a $\nabla$-module relative to $L$.
We call $\calF$ the \emph{generic fibre} of $\calE$;
note that
\begin{equation} \label{eq:compare generic}
T(\calE,(1,\dots,1,\rho)) \leq T(\calF,\rho) \qquad (\rho \in (\epsilon,1)).
\end{equation}
\end{defn}

\begin{prop} \label{P:rel unipotent}
Take $S$ as in Theorem~\ref{T:frob ant}.
Let $\calE$ be a $\nabla$-module on $A_K(S)$
such that $T(\calE,(1,\dots,1,\rho)) = 1$ for all $\rho \in (\epsilon,1)$,
and suppose that the generic fibre of $\calE$ is quasi-unipotent.
Then there exists a positive integer $m$ coprime to $p$ such that
$f_{m I_n}^* \calE$ is unipotent on $A_K([1,1]^{n-1}) \times A_K(\eta,1)$,
in the sense of \cite[\S 3.2]{kedlaya-part1}.
\end{prop}
\begin{proof}
By \eqref{eq:compare generic} and Corollary~\ref{C:close break}, 
we can choose $m$ so that the generic fibre of
$f_{m I_n}^* \calE$ is unipotent. The claim then follows from
\cite[Proposition~3.4.3]{kedlaya-part1}.
\end{proof}

\begin{remark}
Although we have defined a differential highest ramification break, we
have not defined a full set of differential ramification breaks, among which
our highest ramification break is the largest number occurring. 
For the present paper, the highest ramification break is enough; for
the construction of the other breaks, see \cite{kedlaya-swan}.
\end{remark}

\section{Local semistable reduction for monomial valuations}
\label{sec:mono}

We conclude by proving local semistable reduction for monomial valuations.

\setcounter{theorem}{0}

\subsection{Monomial valuations}
\label{subsec:monomial}

\begin{defn}
Let $F$ be a finitely generated field over $k$.
A valuation $v$ on $F$ over $k$ is \emph{monomial}
(in the sense of
\cite[Definition~2.5.3]{kedlaya-part2}) if
\[
\rank(v) = 1, \qquad \ratrank(v) = \trdeg(F/k), \qquad \kappa_v = k.
\]
Note that $v$ is then \emph{minimal} in the sense of
\cite[Definition~4.3.2]{kedlaya-part2}.
Moreover, $v$ is an Abyhankar valuation in the sense of
\cite[Definition~2.5.3]{kedlaya-part2}, which forces the value group of $v$
to be a a finite free $\ZZ$-module.
\end{defn}

\begin{prop}
Let $F$ be a finitely generated field over $k$,
let $v$ be a monomial valuation on $F$ with residue field $k$, 
and let $x_1, \dots, x_n \in F$
be such that $v(x_1), \dots, v(x_n)$ freely generate the value group
of $v$. Then the completion $\widehat{F}$ is isomorphic to the completion 
$k((x_1, \dots, x_n))_v$
of $k(x_1, \dots, x_n)$ under $v$, i.e., the set
of formal sums $\sum_I a_I x^I$ with $a_I \in k$ such that for any $c \in \RR$,
there are only finitely many indices $I$ with $v(x^I) < c$ and $a_I \neq 0$.
\end{prop}
\begin{proof}
(For properties of valuations used in this argument, see for instance
\cite[Chapter~6]{ribenboim}.)
The extension $\widehat{F}$ of $k((x_1, \dots, x_n))_v$ is
finitely generated and of transcendence degree 0, and hence finite.
Suppose this extension is nontrivial.
Since it is immediate (it changes neither the value group 
nor the residue
field), by Ostrowski's theorem 
\cite[Theorem~6.1.2]{ribenboim},
its degree is a power of $p$, as is the degree of its Galois
closure. By an elementary argument with $p$-groups, 
$\widehat{F}$ contains an Artin-Schreier subextension which 
is also immediate.

However, any Artin-Schreier extension of $k((x_1,\dots,x_n))_v$ can
be written as $z^p - z = P(x_1, \dots, x_n)$, where no monomial of
$P$ of negative degree is a $p$-th power. Hence one of the following
is true, yielding a contradiction.
\begin{itemize}
\item
We have $v(P) \geq 0$, in which case the extension is unramified and hence not
immediate.
\item
We have $v(P) < 0$, and 
the lowest degree monomial of $P$ has valuation not divisible by $p$
in the value group; then the extension has strictly larger value group,
so is not immediate.
\item
We have $v(P) < 0$, and 
the lowest degree monomial of $P$ has valuation  divisible by $p$,
but its coefficient is not a $p$-th power in $k$; then the extension
has strictly larger residue field, so is not immediate.
\end{itemize}
This yields the desired result.
\end{proof}

Since monomial valuations are Abyhankar valuations, they satisfy
local uniformization; the following is a special case of
\cite[Theorem~1.1]{knaf-kuhlmann}.
\begin{prop} \label{P:local uniform}
Let $F$ be a finitely generated field over a field $k$, let $v$ be
a monomial valuation on $F$, and let $Z$ be a finite subset of
the valuation ring $\gotho_v$. Then there exists an irreducible
$k$-scheme of finite type $X$ with $k(X) = F$, on which $v$ is centered
at a smooth closed point $x$, and a system of parameters $a_1, \dots, a_n$
of $X$ at $x$ such that each $z \in Z$ can be written as a unit in
$\calO_{X,x}$ times a monomial in the $a_i$.
\end{prop}

\begin{defn}
Let $X$ be a smooth irreducible $k$-variety, and let $v$ be a 
monomial valuation on $k(X)$ centered at a point $x \in X$.
We say a system of parameters $a_1, \dots, a_n$ for $X$ at $x$ is
\emph{descriptive} for $v$ if $v(a_1), \dots, v(a_n)$ generate
$v(k(X)^*)$.
\end{defn}

\begin{prop} \label{P:local uniform2}
Let $(X,D)$ be a smooth pair over an algebraically closed field $k$
with $X$ irreducible,
and let $v$ be a monomial valuation on $k(X)$ over $k$ centered on $X$.
Then there exist a smooth pair $(X',D')$,
a birational (regular) morphism $f: X' \to X$,
a point $x' \in X'$, and a system of parameters $a_1, \dots, a_n$ for $X'$ at $x'$, 
such that:
\begin{itemize}
\item
$f^{-1}(D) \subseteq D'$;
\item
$v$ is centered at $x'$;
\item
$a_1, \dots,a_n$ is descriptive for $v$;
\item
each component of $D'$ is the zero locus of one of the $a_i$.
\end{itemize}
\end{prop}
\begin{proof}
We may as well take $X$ to be affine. Take the set $Z$ to contain:
\begin{itemize}
\item[(a)]
a set of generators of the coordinate ring $k[X]$ as a $k$-algebra;
\item[(b)]
a sequence $t_1, \dots, t_n$ such that $v(t_1), \dots, v(t_n)$ freely generate
$v(k(X)^*)$ as a $\ZZ$-module;
\item[(c)]
some functions which cut out the components of $D$ passing through 
the center of $v$ on $X$.
\end{itemize}
Apply Proposition~\ref{P:local uniform}; if we take $X'$ to be a
sufficiently small open affine neighborhood of the center $x'$ of $v$ on 
the resulting scheme, and take $D'$ to be the zero locus of $a_1, \dots, a_n$,
then $(X',D')$ will form a smooth pair.
By (a), there will be a birational regular map $f: X' \to X$.
By (b), $v(a_1), \dots, v(a_n)$ generate $v(k(X)^*)$
as a $\ZZ$-module. 
By (c), we can force $f^{-1}(D)
\subseteq D'$ by possibly shrinking $X'$. 
This yields the desired result.
\end{proof}

\begin{prop} \label{P:local uniform3}
Let $(X,D)$ be a smooth pair over an algebraically closed field $k$ 
with $X$ irreducible,
let $v$ be a monomial valuation on $k(X)$ centered at a point $x \in D$, 
let $F$ be a finite Galois extension of $k(X)$, and let $w$ be
an extension of $v$ to $F$.
For $(X',D')$ a toroidal blowup of $(U,U \cap D)$ for some
open neighborhood $U$ of $x$ in $X$, 
write $f:  Y' \to X'$ for the normalization of $X'$ in $F$.
Then it is possible to choose $(X',D')$ such that 
$(Y',f^{-1}(D'))$ is a smooth
pair and $w$ is centered on $Y'$.
\end{prop}
\begin{proof}
We may assume without loss of generality that $x$ is the intersection of
all of the components of $D$. Let $y'$ denote the center of $w$ on $Y'$.

Note that the conclusion implies that in a neighborhood of $y'$,
the pullback of $D'$ to $Y'$ as a Cartier divisor is a $\ZZ$-linear
combination of the components of $f^{-1}(D')$.
Consequently, if
$F'$ is an intermediate field between $k(X)$ and $F$, we can prove the claim
by first passing from $k(X)$ to $F'$ and then from $F'$ to $F$: the point is that
in the second step, the toroidal blowup on the middle variety in the tower is
induced by a toroidal blowup on the bottom variety.

We can write $F/k(X)$ as a tower $F/T/U/k(X)$, where $U/k(X)$ is unramified at $v$,
$T/U$ is totally tamely ramified at $v$, and $F/T$ is a $p$-power extension for $p = \charac(k)$
(or the trivial extension if $\charac(k) = 0$). Moreover, by elementary group theory, 
$F/T$ can be written as a tower of $\ZZ/p\ZZ$-extensions. We may thus reduce
to the cases where $F/k(X)$ is unramified, tamely ramified, or an Artin-Schreier
extension.

There is nothing to check in the unramified case. In the tamely ramified case,
the morphism $Y' \to X'$
is toroidal, so $(Y', f^{-1}(D'))$ is automatically toroidal; it thus suffices
to perform toroidal resolution of singularities \cite{kkms} upstairs, as again
we can mimic the toroidal blowups downstairs. In the
Artin-Schreier case, we have $F = k(X)[z]/(z^p - z - h)$ for some $h \in k(X)$
with $v(h) < 0$. By Proposition~\ref{P:local uniform} (or a direct calculation),
we can choose the blowup $(X',D')$ so that at $x' = f(y')$, 
$h^{-1}$ becomes a unit in $\calO_{X',x'}$ times a product of powers of
local parameters of components of $D'$ at $x'$. Then $(Y', f^{-1}(D'))$ is toroidal,
so again toroidal resolution of singularities yields the claim.
\end{proof}

\begin{remark}
Beware that in Proposition~\ref{P:local uniform3}, the morphism 
$Y' \to X'$ is in general not toroidal when $\charac(k) = p > 0$.
This is already true for curves: consider the covering 
\[
\Spec k[x,t]/(t - x^p - x^{p+1}) \to \Spec k[x].
\]
\end{remark}

\subsection{The contagion of unipotence}

\begin{prop} \label{P:contagion}
Let $S$ be the set of $n$-tuples 
$(\rho^{x_1}, \dots, \rho^{x_{n-1}}, \rho)$ for $\rho$ in some interval
$(\epsilon,1)$ and $x = (x_1, \dots, x_{n-1})$ in some rational polyhedral
subset $U$ of $\RR^{n-1}$.
Let $q$ be a power of $p$.
Let $\calE$ be a $\nabla$-module on $A_K(S)$ equipped with an isomorphism
$\sigma^* \calE \cong \calE$ on $A_K(S^{1/q})$ for some map
$\sigma: A_K(S^{1/q}) \to A_K(S)$ obtained by composing the toric map
$f_{qI_{n}}$ with a $q$-power Frobenius lift on $K$.
Suppose $y \in U$ is such that $1, y_1, \dots, y_{n-1}$ are linearly independent
over $\QQ$,
and 
\[
T(\calE,(\rho^{y_1},\dots,\rho^{y_{n-1}},\rho)) = 1 \qquad \mbox{for $\rho
\in (\epsilon,1)$ sufficiently close to $1$}.
\]
Then
there exists a neighborhood $V$ of $y$ in $U$ such that for $x \in V$,
$T(\calE,(\rho^{x_1},\dots,\rho^{x_{n-1}},\rho)) = 1$ for 
$\rho \in (\epsilon,1)$ sufficiently close to $1$.
\end{prop}
\begin{proof}
By Theorem~\ref{T:rational} (applicable because of Remark~\ref{R:Frobenius}),
for each $x \in U \cap \QQ^{n-1}$,
there exists $f(x) \geq 0$ with
\[
f(x) \in \frac{1}{\rank(\calE)!} (\ZZ + x_1 \ZZ + \cdots + x_{n-1} \ZZ).
\]
such that
\[
T(\calE,(\rho^{x_1},\dots,\rho^{x_{n-1}},\rho)) = \rho^{f(x)} 
\qquad \mbox{for $\rho
\in (\epsilon,1)$ sufficiently close to $1$}.
\]
Moreover, $f(x)$ is convex by Proposition~\ref{P:log concave break}.
Thus we may apply Theorem~\ref{T:internally integral}
to deduce that $\rank(\calE)! f$ is internally integral polyhedral.

The boundaries between domains of affinity of $f$ all lie on rational
hyperplanes, whereas $y$ lies on no such hyperplanes because
$1, y_1, \dots, y_{n-1}$ are linearly independent  over $\QQ$. 
Hence $y$ lies in the interior
of some domain of affinity. In that domain, there exist $a_1, \dots, a_{n-1}, 
b \in \ZZ$ such that
\[
\rank(\calE)! f(x) = a_1 x_1 + \cdots + a_{n-1} x_{n-1} + b.
\]
Since $f(y) = 0$ and $1, y_1,\dots,y_{n-1}$ are linearly independent over
$\QQ$, we must have $a_1 = \cdots = a_{n-1} = b = 0$, that is, $f(x) = 0$
identically in an open neighborhood of $y$, as desired.
\end{proof}

\subsection{$F$-isocrystals near a monomial valuation}

We are now ready to prove our first instances of local semistable reduction
at a minimal valuation on a variety of dimension greater than
1. (The theorem
also applies for $X$ of
dimension 1, but in that case
one can simply apply the usual $p$-adic local monodromy theorem for
the same effect.)

\begin{theorem} \label{T:monomial}
Let $X$ be a smooth irreducible $k$-scheme,
let $\overline{X}$ be a partial compactification of $X$,
and let $\calE$ be an $F$-isocrystal on $X$ overconvergent
along $\overline{X} \setminus X$. Then $\calE$ admits
local semistable reduction at any monomial valuation on $k(X)$
centered on $\overline{X}$.
\end{theorem}
\begin{proof}
We may assume $k$ is algebraically closed thanks to 
\cite[Proposition~3.2.6]{kedlaya-part2}. 
Let $v$ be a monomial valuation on $k(X)$.
By Proposition~\ref{P:local uniform2}, there is a smooth pair
$(Y,D)$ containing an open dense subscheme of $X$, such that $v$
is centered at an intersection of components of $D$, and the valuations
of some system of parameters $t_1, \dots, t_n$
at that point freely generate $v(k(X)^*)$.

Put $y_i = v(t_i)/v(t_n)$ for $i=1, \dots, n-1$;
we can then realize $\calE$ as a $\nabla$-module
on $A_K(S)$ for some set $S$ containing $(\rho^{x_1},\dots,\rho^{x_{n-1}},\rho)$ 
for $\rho$ in some interval $(\epsilon,1)$ and 
$x = (x_1, \dots, x_{n-1})$ in some neighborhood of $y$ in $\RR^{n-1}$.
Moreover, $\calE$ admits a Frobenius action for a Frobenius lift on
$A_K(S)$ given by composing a $q$-power
Frobenius lift on $K$ with the toric map $f_{qI_n}$.
Take $\lambda = (y_1,\dots,y_{n-1},1)$ 
and form the $(F, \nabla)$-module $M_\calE$
over $\calR^\lambda$ corresponding to $\calE$.
If $M_\calE$ is unipotent, we may 
apply Proposition~\ref{P:contagion} to deduce that
for $x$ in a possibly smaller neighborhood of $y$,
$T(\calE,(\rho^{x_1}, \dots, \rho^{x_{n-1}},\rho)) = 1$ 
for $\rho$ sufficiently close to 1.

This means (by virtue of Proposition~\ref{P:rel unipotent} applied after a toric
coordinate change)
that by passing to a suitable toroidal blowup in the sense of
\cite{kkms}, we can obtain another smooth pair
$(Y',D')$ such that $v$ is centered at the
intersection of $n$ components of $D'$, and $\calE$ becomes unipotent along
each of those components after making a suitable tamely ramified cover.
(For instance, it suffices to perform a blowup corresponding to a
barycentric subdivision sufficiently many times.)
If we take $m$ sufficiently divisible and prime to $p$, then
pass to a cover that is tamely ramified of degree $m$ along each of the
$n$ components of $D'$, we get a smooth pair $(Y'',D'')$ on which $v$ is centered at an
intersection of components of $D''$, along each of which $\calE$ is unipotent.
By \cite[Theorem~6.4.5]{kedlaya-part1},
$\calE$ extends to a log-isocrystal with nilpotent residues
on $(Y'',D'')$.

If $M_\calE$ is not unipotent,
we apply Theorem~\ref{T:fake annuli local mono} (to produce a good finite cover)
and Proposition~\ref{P:local uniform3} (to toroidalize)
to deduce that after passing up to a suitable quasi-resolution,
we get into the situation where $M_\calE$ is indeed unipotent.
This yields local semistable reduction at $v$, as desired.
\end{proof}

By virtue of earlier work, we obtain the same conclusion more
generally for Abhyankar valuations.
\begin{cor}
Let $X$ be a smooth irreducible $k$-scheme,
let $\overline{X}$ be a partial compactification of $X$,
let $\calE$ be an $F$-isocrystal on $X$ overconvergent
along $\overline{X} \setminus X$, and let $v$ be any Abhyankar 
valuation on $k(X)$ centered on $\overline{X}$.
Then $\calE$ admits local semistable reduction at $v$.
\end{cor}
\begin{proof}
This follows from Theorem~\ref{T:monomial} as in
the proofs of \cite[Proposition~4.2.4 and
Theorem~4.3.4]{kedlaya-part2}.
\end{proof}

\appendix
\section{Some examples}

In this appendix, we make good on two promises of examples
to illustrate aspects of the semistable reduction problem.

\subsection{Finite covers are not enough}

The following example illustrates that one cannot necessarily render unipotent
the local monodromy of an overconvergent $F$-isocrystal by pulling
back along a finite cover instead of an alteration, as alluded to in
the introduction of \cite{kedlaya-part1}.

\begin{example}
Let $\calF$ be the pullback along the map $t \mapsto t^{-1}$ of the
Bessel isocrystal on $\GG_m$, as defined in \cite[Example~6.2.6]{tsuzuki-slope}.
Then there exists a finite flat morphism $f: X \to \PP^1_k$ such that
$f^* \calF$ extends to a convergent log-isocrystal $\calF_1$ on 
$(X, f^{-1}(\{0,\infty\}))$, and the Frobenius slopes of 
$\calF_1$ at a closed point $x \in X$ equal $1/2, 1/2$ if $f(x) = \infty$
and $0,1$ otherwise.

Let $\pi_1,\pi_2: \PP^1_k \times \PP^1_k \to \PP^1_k$ denote the canonical
projections, and put $\calE = \pi_1^* \calF \otimes \pi_2^* \calF$.
Based on the properties of $\calF$, we know that there exists
an alteration $f_1: X_1 \to \PP^1_k \times \PP^1_k$ such that
$f_1^* \calE$ extends to a convergent log-isocrystal $\calE_1$ 
on $X_1$ for some
log structure. Moreover, for one such alteration,
the Frobenius slopes of $\calE_1$ at a closed point $x \in X_1$ equal
\begin{equation} \label{eq:bessel slopes}
\begin{cases} 1,1,1,1 & f_1(x) = (\infty, \infty) \\
1/2, 1/2, 3/2, 3/2 & f_1(x) \in (\{\infty\} \times \AAA^1_k) \cup (\AAA^1_k
\times \{\infty\}) \\
0,1,1,2 & f_1(x) \in \AAA^1_k \times \AAA^1_k;
\end{cases}
\end{equation}
it follows that the same holds for \emph{any} such alteration. (Given a second
such alteration $f_2: X_2 \to \PP^1_k \times \PP^1_k$, we can construct a
third alteration $f_3: X_3 \to \PP^1_k \times \PP^1_k$ factoring through both
$f_1$ and $f_2$, then transfer the information about the Frobenius slopes
from $X_1$ to $X_3$ to $X_2$.)

We now wish to argue that there cannot exist a finite morphism
$f: X \to \PP^2_k$ such that $f^* \calE$ extends to a convergent log-isocrystal
on $X$ for some log structure. To see this, we may reduce to the case
where $f$ is Galois (by replacing the cover by its normal closure), in 
which case the Frobenius slopes of the extension of $f^* \calE$ at a point
$x \in X$ depend only on the projection $f(x)$.

Let $P$ be the closure of the
graph of a rational map $\PP^1_k \times \PP^1_k \dashrightarrow \PP^2_k$
identifying $\AAA^1_k \times \AAA^1_k$ with $\AAA^2_k$. 
Put $Y = X \times_{\PP^2_k} P$, so that base change induces a finite morphism
$f: Y \to P$, and let $f_1$ denote the composition $Y \stackrel{f}{\to} P 
\to \PP^1_k \times \PP^1_k$.
Then the above analysis shows that the Frobenius slopes of the extension of
$f_1^* \calE$ at a point $y \in Y$ depend only on $f_1(y)$.

However, this yields a contradiction as follows. Each of the three components
of $Z = P \setminus \AAA^2_k$ is contracted by one of the projections
$P \to \PP^1_k \times \PP^1_k$ or $P \to \PP^2_k$. Consequently,
the Frobenius slopes must be constant along each component; since $Z$
is connected, the slopes must be constant along all of $f_1^{-1}(Z)$.
However,
this contradicts the explicit formula \eqref{eq:bessel slopes}.
\end{example}

\begin{remark}
This example is not meant to suggest that one is compelled to blow up in the
locus where the isocrystal is already defined. Indeed, it is entirely
possible that one can always use an alteration which is finite \'etale
over that locus; however, even if one had as strong a form of
resolution of singularities in positive characteristic as desired, it 
is not clear how to use the valuation-theoretic approach 
to prove this refined form of semistable reduction.
\end{remark}

\subsection{Extra monodromy on exceptional divisors}

The following example illustrates that one cannot necessarily render unipotent
the local monodromy of an overconvergent $F$-isocrystal by doing so only
for the divisors in a specified good compactification of the locus of
definition, as alluded to in the introduction of this paper.

\begin{example} \label{exa:extra monodromy}
Consider an affine plane $\AAA^2_k$ with coordinates $x,y$, embed it
into a projective plane $\PP^2_k$, and
let $X$ be the complement of the line $y=0$ in $\PP^2_k$.
View $P(x,y,z) = yz^{p^2} - x^{p-1} z^p + z$ as a polynomial in
$k(x,y)[z]$.
One checks that the extension $k(x,y)[z]/(P)$ of
$k(x,y) = k(X)$ defines a finite \'etale cover
$f: Y \to X$. Let $\calE$ be the overconvergent $F$-isocrystal
$f_* \calO_Y$ on $X$. 
We consider twisted polynomials again as in \cite{ore}, but for the
Frobenius automorphism instead of for a derivation.
Over the $y$-adic completion $k(x)((y))$
of $k(x,y)$, we can factor
the twisted polynomial 
$Q = y F^2 - x^{p-1} F + 1$ as $(yF - c)(F - 1/c)$ for 
some $c \equiv x^{p-1} \pmod{y}$; in particular, $c$ has a $(p-1)$-st root
in $k(x)((y))$. We may thus split $P$
over an Artin-Schreier extension of $k(x)((y))$; by Krasner's lemma,
we can realize this as the completion of a degree $p$ extension of $k(x,y)$.

This means that we can construct a finite flat
morphism $g: Y_1 \to \PP^2_k$ of degree $p$ such that $g^* \calE$ has constant
local monodromy along each component of the proper transform of the line
$y=0$. However, if we blow up at $x=y=0$ and complete the function field
along the resulting exceptional divisor, we obtain $k(x/y)((y))$, over which
$Q$ remains irreducible. Consequently, $g^* \calE$ cannot have constant
local monodromy along the proper transform of the exceptional divisor.
\end{example}

\begin{remark}
In Example~\ref{exa:extra monodromy}, the
overconvergent $F$-isocrystal $\calE$ is unit-root because
it is a pushforward of the unit-root isocrystal $\calO_Y$.
Hence one can recover semi\-stable reduction for $\calE$ using
results of Tsuzuki \cite{tsuzuki-duke}. The method of proof follows
the model one would use in the $\ell$-adic setting: convert $\calE$ into
a $p$-adic representation of the \'etale fundamental group of $X$, choose
a stable lattice, and pick a finite \'etale cover of $X$ that trivializes
a suitable quotient of the lattice. Unfortunately, without a unit-root
condition, one has no useful functor from isocrystals to Galois 
representations; the compactness
of the Riemann-Zariski space serves as a replacement for this construction.
\end{remark}


\begin{thebibliography}{99}

\bibitem{amice}
Y. Amice, \textit{Les nombres $p$-adiques}, Presses Universitaires
de France, Paris, 1975.

\bibitem{andre}
Y. Andr\'e, Filtrations de type Hasse-Arf et monodromie $p$-adique,
\textit{Invent. Math.} \textbf{148} (2002), 285--317.

\bibitem{berkovich}
V.G. Berkovich, \textit{Spectral theory and analytic geometry
over non-archimedean fields} (translated by N.I. Koblitz), 
Math. Surveys and Monographs 33, Amer. Math. Soc., Providence, 1990.

\bibitem{christol-dwork}
G. Christol and B. Dwork, Modules diff\'erentiels sur des couronnes,
\textit{Ann. Inst. Fourier (Grenoble)} \textbf{44} (1994), 663--701.

\bibitem{cm3}
G. Christol and Z. Mebkhout, Sur le th\'eor\`eme de l'indice des
\'equations diff\'erentielles $p$-adiques. III,
\textit{Annals of Math.} \textbf{151} (2000), 385--457.

\bibitem{christol-robba}
G. Christol and P. Robba, \textit{\'Equations diff\'erentielles
$p$-adiques, Applications aux sommes exponentielles},
Actualit\'es Math\'ematiques, Hermann, Paris, 1994.

\bibitem{dejong}
A.J. de Jong, Smoothness, semi-stability and alterations,
\textit{Inst. Hautes \'Etudes Sci. Publ. Math.}
\textbf{83} (1996), 51--93.

\bibitem{g-functions}
B. Dwork, G. Gerotto, and F.J. Sullivan, \textit{An introduction
to $G$-functions}, Annals of Math. Studies 133, Princeton Univ. Press,
Princeton, 1994.

\bibitem{dwork-robba}
B. Dwork and P. Robba, On ordinary linear $p$-adic differential equations,
\textit{Trans. Amer. Math. Soc.} \textbf{231} (1977), 1--46.

\bibitem{fujishige}
S. Fujishige, \textit{Submodular functions and optimization}, second
edition, Annals of Discrete Math. 58, Elsevier, Amsterdam, 2005.

\bibitem{kedlaya-local}
K.S. Kedlaya, A $p$-adic local monodromy theorem, \textit{Annals of Math.}
\textbf{160} (2004), 93--184.

\bibitem{kedlaya-mono-over}
K.S. Kedlaya,
Local monodromy for $p$-adic differential equations: an overview,
\textit{Intl. J. of Number Theory} \textbf{1} (2005), 109--154.

\bibitem{kedlaya-slope}
K.S. Kedlaya,
Slope filtrations revisited, 
\textit{Doc. Math.} \textbf{10} (2005), 447--525.

\bibitem{kedlaya-fake}
K.S. Kedlaya, The $p$-adic local monodromy theorem for fake annuli,
\textit{Rend. Sem. Mat. Padova} \textbf{118} (2007), 101--146. 

\bibitem{kedlaya-swan}
K.S. Kedlaya, Swan conductors for $p$-adic differential modules, I:
A local construction, \textit{Alg. and Number Theory} \textbf{1} 
(2007), 269--300. 

\bibitem{kedlaya-part1}
K.S. Kedlaya, Semistable reduction for overconvergent $F$-isocrystals, I:
    Unipotence and logarithmic extensions,
\textit{Compos. Math.} \textbf{143} (2007), 1164--1212.

\bibitem{kedlaya-part2}
K.S. Kedlaya, Semistable reduction for overconvergent $F$-isocrystals, II:
A valuation-theoretic approach, 
\textit{Compos. Math.} \textbf{144} (2008), 657--672. 

\bibitem{kedlaya-rel}
K.S. Kedlaya, Slope filtrations for relative Frobenius, 
\textit{Ast\'erisque}, to appear;
arXiv preprint \texttt{math.NT/0609272v2} (2007).

\bibitem{kedlaya-course}
K.S. Kedlaya, $p$-adic differential equations (version of 8 Feb 08), 
preprint available at
\texttt{http://math.mit.edu/\~{}kedlaya/papers/}.

\bibitem{kkms}
G. Kempf, F.F. Knudsen, D. Mumford, and B. Saint-Donat,
\textit{Toroidal embeddings. I}, Lecture Notes in Math. 339,
Springer-Verlag, Berlin, 1973.

\bibitem{knaf-kuhlmann}
H. Knaf and F.-V. Kuhlmann, Abhyankar places admit local uniformization
in any characteristic, \textit{Ann. Scient. \'Ec. Norm. Sup.}
\textbf{38} (2005), 833--846.

\bibitem{kuhlmann}
F.-V. Kuhlmann, Places of algebraic function fields in arbitrary
characteristic, \textit{Adv. Math.} \textbf{188} (2004), 399--424.

\bibitem{mebkhout}
Z. Mebkhout, Analogue $p$-adique du Th\'eor\`eme de Turrittin
et le Th\'eor\`eme de la monodromie $p$-adique,
\textit{Invent. Math.} \textbf{148} (2002), 319--351.

\bibitem{ore}
O. Ore, Theory of non-commutative polynomials,
\textit{Ann. Math.} \textbf{34} (1933), 480--508.

\bibitem{polya-szego}
G. P\'olya and G. Szeg\H{o}, \textit{Problems and theorems in analysis},
part I, reprint of the 1978 edition, Springer-Verlag, Berlin, 1998.

\bibitem{ribenboim}
P. Ribenboim, \textit{The theory of classical valuations},
Springer-Verlag, New York, 1999.

\bibitem{robba}
P. Robba, Lemmes de Hensel pour les op\'erateurs diff\'erentiels.
Application a la r\'eduction formelle des \'equations diff\'erentielles,
\textit{Enseign. Math.} (2) \textbf{26} (1980), 279--311.

\bibitem{rock}
R.T. Rockafellar, \textit{Convex analysis}, Princeton Univ. Press,
Princeton, 1970.

\bibitem{sabbah}
C. Sabbah,
\'Equations diff\'erentielles \`a points singuliers irr\'eguliers 
et ph\'enom\`ene de Stokes en dimension 2,
\textit{Ast\'erisque} \textbf{263} (2000).

\bibitem{tsuzuki-slope}
N. Tsuzuki,
Slope filtration of quasi-unipotent overconvergent
            $F$-isocrystals,
\textit{Ann. Inst. Fourier (Grenoble)}
\textbf{48} (1998), 379--412.

\bibitem{tsuzuki-duke}
N. Tsuzuki,
Morphisms of $F$-isocrystals and the finite monodromy theorem
            for unit-root $F$-isocrystals,
\textit{Duke Math. J.} \textbf{111} (2002),
385--418.

\bibitem{young-tams}
P.T. Young, Radii of convergence and index for $p$-adic differential
operators, \textit{Trans. Amer. Math. Soc.} \textbf{333} (1992), 769--785.

\end{thebibliography}
\end{document}